\documentclass{scrartcl}
\usepackage[top=3cm,bottom=3cm,left=2.5cm,right=2.5cm]{geometry} 
\usepackage{graphicx} 

\usepackage{amsmath,amssymb,amsthm,amscd} 
\usepackage{empheq} 
\usepackage{subfig} 
\usepackage{here}

\usepackage{xcolor}

\theoremstyle{definition}

\newtheorem{definition}{Definition}

\newtheorem{lemma}{Lemma}

\newtheorem*{remark}{Remark}

\usepackage{hyperref} 

\title{
Numerical computations of split Bregman method for fourth order total variation flow
}
\author{
Yoshikazu Giga\thanks{Graduate School of Mathematical Sciences, The University of
Tokyo, Komaba 3-8-1, Meguro, Tokyo 153-8914, Japan. \textit{E-mail}: \texttt{labgiga@ms.u-tokyo.ac.jp}} \ and
Yuki Ueda\thanks{Department of Applied Mathematics, The Hong Kong Polytechnic University, Hung Hom, Kowloon, Hong Kong. \textit{E-mail}: \texttt{yuki.ueda@polyu.edu.hk}}
}

\begin{document}
\maketitle
\begin{abstract}
The split Bregman framework for Osher-Sol\'e-Vese (OSV) model and fourth order total variation flow are studied. 
We discretize the problem by piecewise constant function and compute $\nabla(-\Delta_{\mathrm{av}})^{-1}$ approximately and exactly. 
Furthermore, we provide a new shrinkage operator for Spohn's fourth order model.
Numerical experiments are demonstrated for fourth order problems under periodic boundary condition.
\end{abstract}


\bigskip

{\noindent \textbf{2010 Mathematics Subject Classification:}
35K25, 35K67, 65M06

\section{Introduction}
A gradient flow has been of great interest in mathematics and mathematical physics because several evolution equations can be regarded as the gradient flows. For example, mathematical models in materials sciences including the Allen-Cahn equation and mean curvature flow can be regarded as second order $L^2$-gradient flows. The Cahn-Hilliard equation is interpreted as a fourth-order $H^{-1}$-gradient flow.

We are interested in several important examples of gradient flows which are of the form
\begin{equation}
\label{Eq:GradientFlow_for_Hilbert}
\dfrac{\partial u}{\partial t} \in -\partial_H E(u)\mbox{ for }t>0,
\end{equation}
where $H$ is a Hilbert space, $E:H\to\mathbb{R}\cup\{\infty\}$ is a convex, lower semi-continuous functional and the subdifferential $\partial_H$ is defined as
\begin{equation}
\partial_HE(u) = \left\{p\in H : E(v)-E(u)\ge(p,v-u)_H\mbox{ for all }v\in H\right\}.
\end{equation} 
In this paper, we consider gradient flows \eqref{Eq:GradientFlow_for_Hilbert} with convex energy $E$ but may be very singular.
We give a few examples. Spohn \cite{Sp93} has proposed a mathematical model for the relaxation of a crystalline surface below the roughening temperature;
\begin{equation}
\label{Eq:SpohnFourthOrder}
u_t = -\Delta\left(\operatorname{div}\left(\beta\dfrac{\nabla u}{|\nabla u|}+|\nabla u|^{p-2}\nabla u\right)\right),
\end{equation}
where $\beta>0$ and $p>1$.
Kashima \cite{Kas04} has presented this model as a fourth order $H^{-1}$-gradient flow for energy functional 
\begin{equation}
\label{Eq:SpohnModel}
E(u) = \beta\displaystyle\int_{\Omega}|Du| + \dfrac{1}{p}\int_{\Omega}|Du|^p.
\end{equation}

Furthermore, the total variation flow, which is the gradient flow for total variation energy, has been studied well in image processing. In 1992, Rudin, Osher and Fatemi \cite{ROF92} have introduced the total variation to image denoising and reconstruction. Their model, which is known as the Rudin-Osher-Fatemi (ROF) model, is described as
\begin{equation}
u = \mathop{\mathrm{argmin}}_{u\in L^2(\Omega)} \left\{\displaystyle\int_{\Omega}|Du|+\dfrac{\lambda}{2}\|u-f\|_{L^2(\Omega)}^2\right\},
\end{equation}
where $\Omega\subset\mathbb{R}^2$ is bounded domain and $f:\Omega\to\mathbb{R}$ is a given noisy image. 
This introduces the second order total variation flow
\begin{equation} 
u_t = \operatorname{div}\left(\dfrac{\nabla u}{|\nabla u|}\right) + \lambda(u-f).
\end{equation}
On the other hand, Osher, Sol\'{e} and Vese \cite{OSV03} have introduced the $H^{-1}$ fidelity and provided Osher-Sol\'e-Vese (OSV) model
\begin{equation}
\label{OSVmodel}
u=\mathop{\mathrm{argmin}}_{u\in H^{-1}(\Omega)}\left\{\displaystyle\int_{\Omega}|Du|+\dfrac{\lambda}{2}\|u-f\|_{H^{-1}(\Omega)}^2\right\},
\end{equation}
where $H^{-1}(\Omega) = (H^1_0(\Omega))^*$. Their model performs better on textured or oscillatory images. Equation \eqref{OSVmodel} gives the fourth order total variation flow
\begin{equation}
\label{Eq:FourthOrder}
u_t = -\Delta\left(\operatorname{div}\left(\dfrac{\nabla u}{|\nabla u|}\right)\right)+\lambda(u-f).
\end{equation}

Performing numerical computations for the ROF model, the OSV model and total variation flow have difficulties because of its singularity. Several studies have suggested numerical schemes for the ROF model and second order total variation flow. Especially, the split Bregman framework is well-known as an efficient solver for the ROF model.
The aim of this paper is to provide a new numerical scheme, which is based on the backward Euler method and split Bregman framework, for fourth order total variation flow and Spohn's fourth order model. Numerical experiments are demonstrated for fourth order problems under periodic boundary condition. 

The split Bregman method, which is based on the Bregman iterative scheme \cite{Br'e67}, has been studied and performed in image processing (for example, see \cite{OBG05}).
Goldstein and Osher \cite{GO09} have proposed the alternating split Bregman method. Their method separates the ``$L^1$'' minimization part and ``$L^2$'' part. 
The alternating split Bregman method has several advantages. They have mentioned that the ``$L^2$'' part is differentiable, and the shrinking method can be applied to the ``$L^1$'' part for the ROF model. Therefore it is extremely efficient solver and easy to code. The split Bregman framework can be performed for second order total variation flow easily. For example, 
\begin{equation} 
\label{Eq:Second_TV_Flow}
u_t = \operatorname{div}\left(\dfrac{\nabla u}{|\nabla u|}\right)
\end{equation}
introduces the subdifferential formulation $u_t\in -\partial F(u)$, where $F(u) = \int_{\Omega}|Du|$. 
We let $u_t\approx (u^{k+1}-u^k)/\tau$ and apply the backward Euler method to the subdifferential formulation, then we obtain
\begin{equation}
u^{k+1} = \mathop{\mathrm{argmin}}_{u\in L^2(\Omega)}\left\{\displaystyle\int_{\Omega}|Du|+\dfrac{1}{2\tau}\|u-u^k\|_{L^2(\Omega)}^2\right\},
\end{equation}
where $\tau$ is time step size. This is essentially the same problem as the ROF model, therefore the split Bregman framework can be applied to second order total variation flow.

In this paper, we propose the split Bregman framework for the OSV model \eqref{OSVmodel}, fourth order total variation flow
\begin{equation}
\label{Eq:FourthOrderTVFlow}
u_t = -\Delta\left(\operatorname{div}\left(\dfrac{\nabla u}{|\nabla u|}\right)\right)
\end{equation}
and Spohn's fourth order model \eqref{Eq:SpohnFourthOrder}. For simplicity, we consider one-dimensional torus $\mathbb{T}$. We introduce spatial discretization by piecewise constant functions, then we compute $\nabla(-\Delta_{\mathrm{av}})^{-1}v_h$ for piecewise constant function $v_h$ approximately or exactly. 
We apply the discrete gradient and discrete inverse Laplacian in our first scheme. 
In our second scheme, we calculate the inverse Laplacian for piecewise constant functions directly by using the second degree B-spline, which is continuously differentiable piecewise polynomials. 
The problem can be reduced into a minimization problem on the Euclidean space, which is included in earlier studies for the ROF model.
Therefore we can apply the split Bregman framework to fourth order problems.
Several theoretical results such as the convergence \cite{COS10} can be applied to our scheme directly.
Both of our two schemes are demonstrated for fourth order problems, and we can check that they perform quite well.
Furthermore, we introduce a new shrinkage operator for Spohn's fourth order model. This enables to perform the numerical experiment for the relaxation of a crystalline surface below the roughening temperature quite effectively. Our scheme can be extended to fourth order problems on the two-dimensional torus. We also suggest a shrinkage operator for two-dimensional Spohn's model.

Let us quickly overview some earlier results.
There are many mathematical studies for the second and fourth order total variation flow. The well-posedness for fourth order total variation flow can be proved by considering the right hand side in \eqref{Eq:FourthOrderTVFlow} as a subdifferential of a total variation in $H^{-1}(\Omega)$ (see \cite{Kas04}). This enables us to use the theory of maximal monotone operators \cite{Kom67,Bre73}.
On the other hand, Elliott and Smitheman \cite{ES07} have proved the well-posedness for fourth order total variation flow by using the Galerkin method. 
Adjusting the methods in \cite{GK11}, Giga, Kuroda and Matsuoka \cite{GKM14} have established the extinction time estimate under Dirichlet boundary condition.
Numerical computations which include anisotropic diffusion are performed in \cite{MMR15} for second order models.
Note that even for the second order total variation flow \eqref{Eq:Second_TV_Flow}, because of singularity at $\nabla u=0$, the speed of evolution is determined by nonlocal quantities. Therefore the definition of the solution itself is nontrivial. 
For the second order model, the comparison principle holds, and the theory of viscosity solutions is applicable to show well-posedness for a wide class of total variation type equations \cite{GGP13,GGP14}. However, for the fourth order problem, the comparison principle does not hold in general (see \cite[Theorem 3.7]{GG10}), and the theory of viscosity solutions is not applicable.
For more details of mathematical analysis, we refer the reader to \cite{GG10} and references therein. 

Several studies have considered the fourth order problem under periodic boundary condition.
Kashima \cite{Kas12} has studied the characterization of subdifferential in $H_{\mathrm{av}}^{-1}(\mathbb{T}^d)$.
The exact profile of the fourth order total variation flow has been studied in \cite{GG10}. The extinction time estimate under periodic boundary condition has been established in \cite{GK11}.
A duality based numerical scheme which applies the forward-backward splitting has been proposed in \cite{GMR19}. 
Kohn and Versieux \cite{KV10} have performed the numerical computation for Spohn's model. Their numerical scheme is based on the backward Euler method, mixed finite element method and regularization for singularity. They have proved the convergence by combining the regularization error estimate with temporal and spatial semidiscretization error estimates.

The application of the split Bregman framework to crystalline flow has also been studied through what is called a level-set method.
A crystalline mean curvature flow has been proposed independently in \cite{AG89} and \cite{Tay91}. 
Oberman, Osher, Takei and Tsai \cite{OOTT11} have proposed applying the split Bregman framework to the level-set equation of mean curvature flow. Po\v{z}\'{a}r \cite{Poza18} has studied self-similar solutions of three dimensional crystalline mean curvature flow and presented a numerical scheme which is based on the finite element method and split Bregman framework. However, all calculations given there are for the second order model.

A level-set method for mean curvature flow interprets the motion of curvature flow by a level-set of a solution of 
\begin{equation}
\label{Eq:LevelSetMeanCurvature}
\dfrac{\partial u}{\partial t}-|\nabla u|\operatorname{div}\left(\dfrac{\nabla u}{|\nabla u|}\right)=0.
\end{equation} 
It is a powerful tool to calculate evolution which experiences topological change. It was introduced by Osher and Sathian \cite{OS88} as a numerical way to calculate the mean curvature flow. Note that the level-set mean curvature equation \eqref{Eq:LevelSetMeanCurvature} looks similar to \eqref{Eq:Second_TV_Flow}. However, the singularity of \eqref{Eq:LevelSetMeanCurvature} at $\nabla u=0$ is weaker than one of \eqref{Eq:Second_TV_Flow} because of the multiplier $|\nabla u|$. Therefore it is not necessary to study nonlocal quantities for the level-set mean curvature equation. Its analytic foundation such as well-posedness and comparison principle has been established in \cite{CGG91,ES91}. For more details, we refer the readers to \cite{Gig06}. Very recently, the analytic foundation of the level-set method is extended to crystalline flow by Po\v{z}\'{a}r and the first author \cite{GP16,GP18} and Chambolle, Morini and Ponsiglione \cite{CMP17} and with Novaga \cite{CMNP17} independently. Their methods are quite different.

This paper is organized as follows. Section \ref{Sec:Preliminary} states the definition of $H^{-1}_{\mathrm{av}}(\mathbb{T})$ and the total variation. We introduce the discretization by piecewise constant functions in Section \ref{Sec:Discrete}. Furthermore, we propose two schemes for reducing $\|\cdot\|_{H^{-1}_{\mathrm{av}}(\mathbb{T})}$ into Euclidean norm. Section \ref{Sec:SplitBregman} presents the split Bregman framework for the OSV model and fourth order total variation flow problem. In Section \ref{Sec:Spohn}, we describe the shrinking method for Spohn's model. This report presents numerical examples on one-dimensional torus in Section \ref{Sec:NumExample}. Finally, we extend our scheme to two-dimensional fourth order problems under periodic boundary condition in Section \ref{Sec:TwoDim}.

\section{Preliminary}
\label{Sec:Preliminary}
\subsection{Fourier analysis on the torus $\mathbb{T}$}
First, we review some of the standard facts on the Fourier analysis for one-dimensional torus $\mathbb{T} = \mathbb{R}/\mathbb{Z}$. In this paper, we regard $\mathbb{T}$ as an interval $[0,1]$ with periodic boundary condition. The Fourier transform for $f\in L^1(\mathbb{T})$ and definition of Sobolev space on $\mathbb{R}$ are explained in \cite[Chapter 3]{Gra14a} and \cite[Chapter 1.3]{Gra14b}, respectively. 
Let $\mathcal{D}(\mathbb{T})$ be the complex-valued function space $C^{\infty}(\mathbb{T})$ endowed with the usual test function topology and $\mathcal{D}'(\mathbb{T})$ be its dual. The Fourier coefficient $\widehat{f}_T(\xi)\in \mathbb{C}$ for $f\in \mathcal{D}'(\mathbb{T})$ is defined by the generalized Fourier transform (for example, see \cite[Chapter 5]{Gru09});
\begin{equation}
\widehat{f}_T(\xi) = \langle f,e^{-2\pi i\xi x}\rangle_{\mathcal{D}'(\mathbb{T}),\mathcal{D}(\mathbb{T})}.
\end{equation}
The generalized Fourier transform satisfies similar properties to Fourier transform for $f\in L^1(\mathbb{T})$, for example,
\begin{equation}
\label{Eq:FourierDeriv}
\widehat{df/dx}_T(\xi) = \left\langle f,\dfrac{d}{dx}e^{-2\pi i\xi x}\right\rangle_{\mathcal{D}'(\mathbb{T}),\mathcal{D}(\mathbb{T})} = -2\pi i\xi\widehat{f}_T(\xi)
\end{equation}
for all $f\in\mathcal{D}'(\mathbb{T})$.
Furthermore, the Fourier coefficients $\widehat{f}_T(\xi)\in\mathbb{C}$ satisfies
\begin{equation}
f(x) = \displaystyle\sum_{\xi\in\mathbb{Z}}\widehat{f}_T(\xi)e^{2\pi i\xi x}
\end{equation}
for all $f\in \mathcal{D}'(\mathbb{T})$ (see \cite[Chapter 8.2]{Gru09}).
In this Fourier series, the convergence should be understood in the natural topology of $\mathcal{D}'(\mathbb{T})$.
It is well-known that $\mathcal{D}(\mathbb{T})$ is dense in $L^2(\mathbb{T})$, therefore we have $L^2(\mathbb{T})\simeq (L^2(\mathbb{T}))^*\subset \mathcal{D}'(\mathbb{T})$, where $(L^2(\mathbb{T}))^*$ is the dual space of $L^2(\mathbb{T})$ (for example, see \cite[Chapter 5.2]{Bre11}).
Using the generalized Fourier transform, the Lebesgue space and the Sobolev space on $\mathbb{T}$ are defined as follows:
\begin{equation}
\label{Def:L2T}
L^2(\mathbb{T}) = \left\{f\in \mathcal{D}'(\mathbb{T}) : \displaystyle\sum_{\xi\in\mathbb{Z}}|\widehat{f}_T(\xi)|^2<\infty\right\},
\end{equation}
\begin{equation}
\label{Def:H1T}
H^1(\mathbb{T}) = \left\{f\in L^2(\mathbb{T}) : \displaystyle\sum_{\xi\in\mathbb{Z}}\xi^2|\widehat{f}_T(\xi)|^2<\infty\right\}=\left\{f\in \mathcal{D}'(\mathbb{T}) : \displaystyle\sum_{\xi\in\mathbb{Z}}(1+\xi^2)|\widehat{f}_T(\xi)|^2<\infty\right\},
\end{equation}
\begin{equation}
\label{Def:H-1T}
H^{-1}(\mathbb{T}) = \left\{f\in\mathcal{D}'(\mathbb{T}) : \displaystyle\sum_{\xi\in\mathbb{Z}}(1+\xi^2)^{-1}|\widehat{f}_T(\xi)|^2<\infty\right\}.
\end{equation}
Note that the duality pairing can be described formally as
\begin{equation}
\langle f,g\rangle_{H^{-1}(\mathbb{T}),H^1(\mathbb{T})} = \displaystyle\sum_{\xi\in\mathbb{Z}}\widehat{f}_T(\xi)\overline{\widehat{g}_T(\xi)}  = \displaystyle\int_{\mathbb{T}}f(x)\overline{g(x)}~dx
\end{equation}
for all $f\in H^{-1}(\mathbb{T})$ and $g\in H^1(\mathbb{T})$.

\subsection{The inverse Laplacian $(-\Delta_{\mathrm{av}})^{-1}$}
We consider the functions on $\mathbb{T}$ whose average are equal to zero. Let
\begin{subequations}
\begin{align}
L^2_{\mathrm{av}}(\mathbb{T}) &= \left\{f\in L^2(\mathbb{T}) : \displaystyle\int_{\mathbb{T}}f(x)~dx=0\right\},\\
H^1_{\mathrm{av}}(\mathbb{T}) &= L^2_{\mathrm{av}}(\mathbb{T})\cap H^1(\mathbb{T}) = \left\{f\in H^1(\mathbb{T}) : \displaystyle\int_{\mathbb{T}}f(x)~dx = 0\right\},\\
H^{-1}_{\mathrm{av}}(\mathbb{T}) &= \left\{f\in H^{-1}(\mathbb{T}) : \langle f,1\rangle_{H^{-1}(\mathbb{T}),H^1(\mathbb{T})} =0\right\}.
\end{align}
\end{subequations}
These definitions agree with the following ones:  
\begin{equation}
L^2_{\mathrm{av}}(\mathbb{T}) = \left\{f\in \mathcal{D}'(\mathbb{T}) : \displaystyle\sum_{\xi\neq0}|\widehat{f}_T(\xi)|^2<\infty\mbox{ and }\widehat{f}_T(0)=0\right\},
\end{equation}
\begin{equation}
H^1_{\mathrm{av}}(\mathbb{T}) = \left\{f\in\mathcal{D}'(\mathbb{T}) : \displaystyle\sum_{\xi\neq0}\xi^2|\widehat{f}_T(\xi)|^2<\infty\mbox{ and }\widehat{f}_T(0) = 0\right\},
\end{equation}
\begin{equation}
H^{-1}_{\mathrm{av}}(\mathbb{T}) = \left\{f\in\mathcal{D}'(\mathbb{T}) : \displaystyle\sum_{\xi\neq0}\xi^{-2}|\widehat{f}_T(\xi)|^2<\infty\mbox{ and }\widehat{f}_T(0) = 0\right\}.
\end{equation}
It is easy to check that each of these spaces are Hilbert space endowed with the inner products
\begin{subequations}
\begin{align}
(f,g)_{L^2_{\mathrm{av}}(\mathbb{T})}&=\displaystyle\sum_{\xi\neq0}\widehat{f}_T(\xi)\overline{\widehat{g}_T(\xi)},\\
(f,g)_{H^1_{\mathrm{av}}(\mathbb{T})}&=\displaystyle\sum_{\xi\neq0}4\pi^2\xi^2\widehat{f}_T(\xi)\overline{\widehat{g}_T(\xi)},\\
(f,g)_{H^{-1}_{\mathrm{av}}(\mathbb{T})}&=\displaystyle\sum_{\xi\neq0}\dfrac{1}{4\pi^2}\xi^{-2}\widehat{f}_T(\xi)\overline{\widehat{g}_T(\xi)},
\end{align}
\end{subequations}
respectively. These inner products introduce the norms $\|\cdot\|_{L^2_{\mathrm{av}}(\mathbb{T})}$, $\|\cdot\|_{H^1_{\mathrm{av}}(\mathbb{T})}$ and $\|\cdot\|_{H^{-1}_{\mathrm{av}}(\mathbb{T})}$. It is easy to check that 
\begin{alignat}{2}
\|f\|_{L^2_{\mathrm{av}}(\mathbb{T})}&=\|f\|_{L^2(\mathbb{T})} &\qquad\mbox{for all }f\in L^2_{\mathrm{av}}(\mathbb{T}),\\
\|f\|_{H^1_{\mathrm{av}}(\mathbb{T})} &= \|df/dx\|_{L^2(\mathbb{T})} &\qquad\mbox{for all }f\in H^1_{\mathrm{av}}(\mathbb{T}).
\end{alignat}

Fix $u\in H^1_{\mathrm{av}}(\mathbb{T})$ arbitrarily. Let $c(\xi) = 4\pi^2\xi^2\widehat{u}_T(\xi)\in\mathbb{C}$ and $f(x) = \sum_{\xi\in\mathbb{Z}}c(\xi)e^{2\pi i\xi x}$, then we have $c(0)=0$ and $\sum_{\xi\neq0}\xi^{-2}|c(\xi)|^2 = 16\pi^4\sum_{\xi\neq0}\xi^2|\widehat{u}_T(\xi)|^2<\infty$. This implies $f\in H^{-1}_{\mathrm{av}}(\mathbb{T})$. Moreover, 
\begin{equation}
-\Delta u(x) = \displaystyle\sum_{\xi\in\mathbb{Z}}\widehat{(-\Delta u)}_T(\xi)e^{2\pi i\xi x}
=\displaystyle\sum_{\xi\in\mathbb{Z}}4\pi^2\xi^2\widehat{u}_T(\xi)e^{2\pi i\xi x} = f(x),
\end{equation} 
where $-\Delta u = -d^2u/dx^2$. Consequently, $u=(-\Delta_{\mathrm{av}})^{-1}f$ defines $(-\Delta_{\mathrm{av}})^{-1}:H^1_{\mathrm{av}}(\mathbb{T})\to H^{-1}_{\mathrm{av}}(\mathbb{T})$. We call this operator the \textit{inverse Laplacian}. 
Let $f\in H^{-1}_{\mathrm{av}}(\mathbb{T})$ and $u=(-\Delta)_{\mathrm{av}}^{-1}f\in H^1_{\mathrm{av}}(\mathbb{T})$, then we have
\begin{align*}
\|f\|_{H^{-1}_{\mathrm{av}}(\mathbb{T})}^2 &= \|-\Delta u\|_{H^{-1}_{\mathrm{av}}(\mathbb{T})}^2\\
&= \displaystyle\sum_{\xi\neq0}\dfrac{1}{4\pi^2}\xi^{-2}\widehat{(-\Delta u)}_T(\xi)\overline{\widehat{(-\Delta u)}_T(\xi)}\\
&= \displaystyle\sum_{\xi\neq0}4\pi^2\xi^2|\widehat{u}_T(\xi)|^2 = \|u\|_{H^1_{\mathrm{av}}(\mathbb{T})}^2.
\end{align*}
This implies
\begin{equation}
\label{Lem:H-1av_norm}
\|f\|_{H^{-1}_{\mathrm{av}}(\mathbb{T})} =\|(-\Delta_{\mathrm{av}})^{-1}f\|_{H^1_{\mathrm{av}}(\mathbb{T})}= \|\nabla (-\Delta_{\mathrm{av}})^{-1}f\|_{L^2(\mathbb{T})}
\end{equation}
for all $f\in H^{-1}_{\mathrm{av}}(\mathbb{T})$, where $\nabla = d/dx$.

\subsection{Bounded variation and $H^{-1}$ fidelity for the torus $\mathbb{T}$}
We recall the spaces of functions of \textit{bounded variation} in one-dimensional torus.
\begin{definition}[Definition 3.3.13 of \cite{Gra14a}]
For a measurable function $f$ on $\mathbb{T}$ which is defined everywhere, we define the \textit{total variation} as
\begin{equation}
\displaystyle\int_{\mathbb{T}}|Df| = \operatorname{ess} \sup\left\{\displaystyle\sum_{j=1}^M|f(x_j)-f(x_{j-1})| : 0=x_0<x_1<\dots<x_M=1\right\},
\end{equation}
where the supremum is taken over all partition of the interval $[0,1]$. 
we say $f$ is \textit{bounded variation} if the total variation of $f$ is bounded.
Moreover, we define 
\begin{equation}
BV(\mathbb{T}) = \left\{v\in \mathcal{D}'(\mathbb{T}) : \displaystyle\int_{\mathbb{T}}|Df|<\infty\right\}.
\end{equation}
\end{definition}

\begin{remark}
In the definition, $D$ can be regarded as the distributional derivative, and $Df$ can be identified with a signed Borel measure.
\end{remark}
\begin{remark}
The total variation on a general open set $\Omega\subset \mathbb{R}^d$ is defined as
\begin{equation}
\displaystyle\int_{\Omega}|Dv| = \sup\left\{-\int_{\Omega}u(x)\operatorname{div}\phi(x)~dx : \phi\in C^{\infty}_0(\Omega;\mathbb{R}^d)\mbox{ and }\|\phi\|_{L^{\infty}(\Omega)}\le1\right\},
\end{equation}
and the space of bounded variation is defined as
\begin{equation}
BV(\Omega) = \left\{v\in L^1(\Omega) : \displaystyle\int_{\Omega}|Dv|<\infty\right\}.
\end{equation}
It is well-known that if $v\in W^{1,1}(\Omega)$, then
\begin{equation}
\displaystyle\int_{\Omega}|Dv| = \int_{\Omega}|\nabla v|~dx = |v|_{W^{1,1}(\Omega)},
\end{equation}
and therefore $W^{1,1}(\Omega)\subset BV(\Omega)\subset L^1(\Omega)$.
\end{remark}

We define the functional $\Phi:H^{-1}_{\mathrm{av}}(\mathbb{T})\to \mathbb{R}\cup\{\infty\}$ as
\begin{equation}
\Phi(v) =\left\{\begin{array}{ll}
\displaystyle\int_{\mathbb{T}}|Dv|&\mbox{if }v\in BV(\mathbb{T})\cap H^{-1}_{\mathrm{av}}(\mathbb{T}),\\
\infty&\mbox{otherwise.}\end{array}\right.
\end{equation}
Note that $\Phi:H^{-1}_{\mathrm{av}}(\mathbb{T})\to \mathbb{R}\cup\{\infty\}$ is nonnegative, proper, lower semi-continuous and convex. In this paper, we consider the gradient flow equation of the form
\begin{equation}
\label{Eq:GradientFlow}
(\mbox{gradient flow})\left\{\begin{array}{rll}
\dfrac{du}{dt}(t) &\in-\partial_{H^{-1}_{\mathrm{av}}(\mathbb{T})}\Phi(u(t))&\mbox{for a.e. }t>0,\\
u(\cdot,0) &=u_0\in H^{-1}_{\mathrm{av}}(\mathbb{T}),&
\end{array}\right.
\end{equation}
where the subdifferential $\partial_{H^{-1}_{\mathrm{av}}(\mathbb{T})}$ is defined as
\begin{equation}
\partial_{H^{-1}_{\mathrm{av}}(\mathbb{T})}F(u) = \left\{p\in H^{-1}_{\mathrm{av}}(\mathbb{T}) : F(v)-F(u)\ge (p,v-u)_{H^{-1}_{\mathrm{av}}(\mathbb{T})} \mbox{ for all }v\in H^{-1}_{\mathrm{av}}(\mathbb{T})\right\}
\end{equation}
for any convex functional $F:H^{-1}_{\mathrm{av}}(\mathbb{T})\to \mathbb{R}\cup\{\infty\}$ and $u\in H^{-1}_{\mathrm{av}}(\mathbb{T})$. It is well-known that the theory of maximal monotone operators shows the existence and uniqueness of solution $u\in C([0,\infty),H^{-1}_{\mathrm{av}}(\mathbb{T}))$ to equation \eqref{Eq:GradientFlow} (for example, see \cite{Kom67}).

Let $\tau>0$ be the temporal step size. We consider the backward Euler method for gradient flow equation \eqref{Eq:GradientFlow}; for given $u^k\in H^{-1}_{\mathrm{av}}(\mathbb{T})$, find $u^{k+1}\in H^{-1}_{\mathrm{av}}(\mathbb{T})$ such that
\begin{equation}
\dfrac{u^{k+1}-u^k}{\tau} 
\in -\partial_{H^{-1}_{\mathrm{av}}(\mathbb{T})}\Phi(u^{k+1}).
\end{equation}
This can be reduced to solving the following minimization problem:
\begin{equation}
\label{Eq:BackEuler}
u^{k+1} = \displaystyle\mathop{\mathrm{argmin}}_{u\in H^{-1}_{\mathrm{av}}(\mathbb{T})}\left\{ \Phi(u)+\dfrac{1}{2\tau}\|u-u^k\|_{H^{-1}_{\mathrm{av}}(\mathbb{T})}^2\right\}.
\end{equation}
Since $\Phi$ is convex, such $u^{k+1}$ is uniquely determined.
The convergence of backward Euler method has been proved in \cite{Kom67}.
Note that equation \eqref{Eq:BackEuler} is similar to the OSV model \cite{OSV03} which can be described as
\begin{equation}
\label{Problem:OSV}
(\mbox{OSV})\left\{\begin{array}{l}
\mbox{Find $u\in H^{-1}_{\mathrm{av}}(\mathbb{T})$ such that}\\
u = \mathop{\mathrm{argmin}}\left\{\Phi(u) + \dfrac{\lambda}{2}\|u-f\|_{H^{-1}_{\mathrm{av}}(\mathbb{T})}^2\right\},
\end{array}\right.
\end{equation}
where $f\in H^{-1}_{\mathrm{av}}(\mathbb{T})$ is given data and $\lambda>0$ is an artificial parameter. 
The existence result in convex analysis (for example, see \cite[Cor3.23]{Bre11}) gives that the minimizer $u\in BV(\mathbb{T})\cap H^{-1}_{\mathrm{av}}(\mathbb{T})$ exists. 

Hereafter, we consider the following minimization problem: find $u\in H^{-1}_{\mathrm{av}}(\mathbb{T})$ such that
\begin{equation}
(\mathrm{P}0)\quad \displaystyle\mathop{\mathrm{minimize}}_{u\in H^{-1}_{\mathrm{av}}(\mathbb{T})}\left\{\Phi(u)+\dfrac{\lambda}{2}\|u-f\|_{H^{-1}_{\mathrm{av}}(\mathbb{T})}^2\right\},
\end{equation}
where $f\in H^{-1}_{\mathrm{av}}(\mathbb{T})$ is a given data or $f=u^k$, and $\lambda$ is a given parameter or $\lambda = 1/\tau$. 
This involves both of (OSV) and the backward Euler method for (gradient flow).
Furthermore, $(\mbox{P}0)$ introduces the following constrained problem:
\begin{equation}
\label{Eq:Constrained}
(\mathrm{P}1)\quad \displaystyle\mathop{\mathrm{minimize}}_{u\in H^{-1}_{\mathrm{av}}(\mathbb{T})}\left\{\displaystyle\int_{\mathbb{T}}|d|+\dfrac{\lambda}{2}\|u-f\|_{H^{-1}_{\mathrm{av}}(\mathbb{T})}^2 : d=Du\right\}.
\end{equation}

\begin{remark}
When we consider the Spohn's model
\begin{equation}
u_t = -\Delta\left(\operatorname{div}\left(\beta\dfrac{\nabla u}{|\nabla u|}+|\nabla u|^{p-2}\nabla u\right)\right),
\end{equation}
the subdifferential formulation is given as $u_t \in -\partial_{H^{-1}_{\mathrm{av}}(\mathbb{T})}\widetilde{\Phi}(u)$, where
\begin{equation}
\widetilde{\Phi}(u) = \beta\displaystyle\int_{\mathbb{T}}|Du| + \dfrac{1}{p}\int_{\mathbb{T}}|Du|^p.
\end{equation}
Therefore the backward Euler method yields
\begin{equation}
u^{k+1} = \mathop{\mathrm{argmin}}_{u\in H^{-1}_{\mathrm{av}}(\mathbb{T})}\left\{\widetilde{\Phi}(u)+\dfrac{1}{2\tau}\|u-u^k\|_{H^{-1}_{\mathrm{av}}(\mathbb{T})}^2\right\}.
\end{equation}
Then we consider the constraint problem 
\begin{equation}
\mathop{\mathrm{minimize}}_{u\in H^{-1}_{\mathrm{av}}(\mathbb{T})}\left\{\beta\displaystyle\int_{\mathbb{T}}|d|+\dfrac{1}{p}\int_{\mathbb{T}}|d|^p + \dfrac{\lambda}{2}\|u-f\|_{H^{-1}_{\mathrm{av}}(\mathbb{T})}^2 : d=Du\right\}.
\end{equation}
\end{remark}

\section{Discretization for total variation flow and OSV model} 
\label{Sec:Discrete}
\subsection{Discretization for minimization problem}
We introduce the (spatial) discretization for the problem $(\mbox{P}1)$. Let $N\in\mathbb{N}$ be the partition number, $h=1/N$ and $x_n = nh$. We regard $x_0 = x_N$, then $\{x_n\}_{n=0}^N$ gives an uniform partition for $\mathbb{T}$. Furthermore, we let $x_{n+1/2} = (x_n+x_{n+1})/2 = (n+1/2)h$ for $n = -1,0,\dots, N$, where $x_{-1/2}$ and $x_{N+1/2}$ are identified with $x_{N-1/2}$ and $x_{1/2}$, respectively. Then we define the following spaces of piecewise constant functions:
\begin{subequations}
\begin{align}
V_{h} &= \left\{v_h:\mathbb{T}\to\mathbb{R} : v_h|_{I_n}\in \mathbb{P}_0(I_n)\mbox{ for all }n=0,\dots,N\right\},\\
V_{h0} &= \left\{v_h = \displaystyle\sum_{n=1}^Nv_n\boldsymbol{1}_{I_n}\in V_h : \displaystyle\sum_{n=1}^Nv_n = 0\right\},\\
\widehat{V}_{h} &= \{d_h:I\to\mathbb{R} : d_h|_{[x_{n-1},x_n)}\in\mathbb{P}_0([x_{n-1},x_n))\mbox{ for all }n=1,\dots,N\},
\end{align}
\end{subequations}
where $I=[0,1]$, $I_n = [x_{n-1/2},x_{n+1/2})$, $\mathbb{P}_0(I_n)$ is a space of constant functions on interval $I_n$ and $\boldsymbol{1}_{I_{n}}$ is its characteristic function. 
Note that $V_{h0}$ is a finite dimensional subspace of $H^{-1}_{\mathrm{av}}(\mathbb{T})$.

Furthermore, we define $D_h:V_{h0}\to\widehat{V}_h\cap L^2_{\mathrm{av}}(I)$ as
\begin{equation}
D_hv_h = \displaystyle\sum_{n=1}^N(v_{n}-v_{n-1})\boldsymbol{1}_{[x_{n-1},x_n)},
\end{equation}
where $v_0$ is identified with $v_N$.
Let $d_h=D_hv_h\in\widehat{V}_h$, $\textbf{d} = (d_1,\dots, d_N)^{\mathrm{T}}\in\mathbb{R}^N$ for $d_h=\sum_{n=1}^Nd_n\boldsymbol{1}_{[x_{n-1},x_n)}$, $\widetilde{\textbf{v}} = (v_1,v_2,\dots, v_{N})^{\mathrm{T}}\in\mathbb{R}^N$ for $v_h=\sum_{n=1}^Nv_n\boldsymbol{1}_{I_n}\in V_{h0}$, then we have
\begin{equation}
\Phi(v_h) = \dfrac{\|D_hv_h\|_{L^1(I)}}{h} = \dfrac{\|d_h\|_{L^1(I)}}{h} = h\|\nabla_h\widetilde{\textbf{v}}\|_1 = \|\textbf{d}\|_1,
\end{equation}
where $\nabla_h:\mathbb{R}^N\to\mathbb{R}^N$ is the discrete gradient
\begin{equation}
\nabla_h = h^{-1}\begin{pmatrix}
1&0&\dots&0&-1\\
-1&1&\dots&0&0\\
\vdots&&\ddots&&\vdots\\
0&0&\dots&-1&1
.\end{pmatrix}\in\mathbb{R}^{N\times N}.
\end{equation}
Note that $D_hv_h\in \widehat{V}_h\cap L^2_{\mathrm{av}}(I)\subset L^2_{\mathrm{av}}(I)$ for all $v_h\in V_{h0}$; however, $D_hv_h\not\in H^{-1}_{\mathrm{av}}(\mathbb{T})$ because it does not satisfy the periodic boundary condition (see Figure \ref{Fig:vh_and_Dhvh}). 
\begin{figure}[tb]
\centering
\includegraphics[clip,scale=.7]{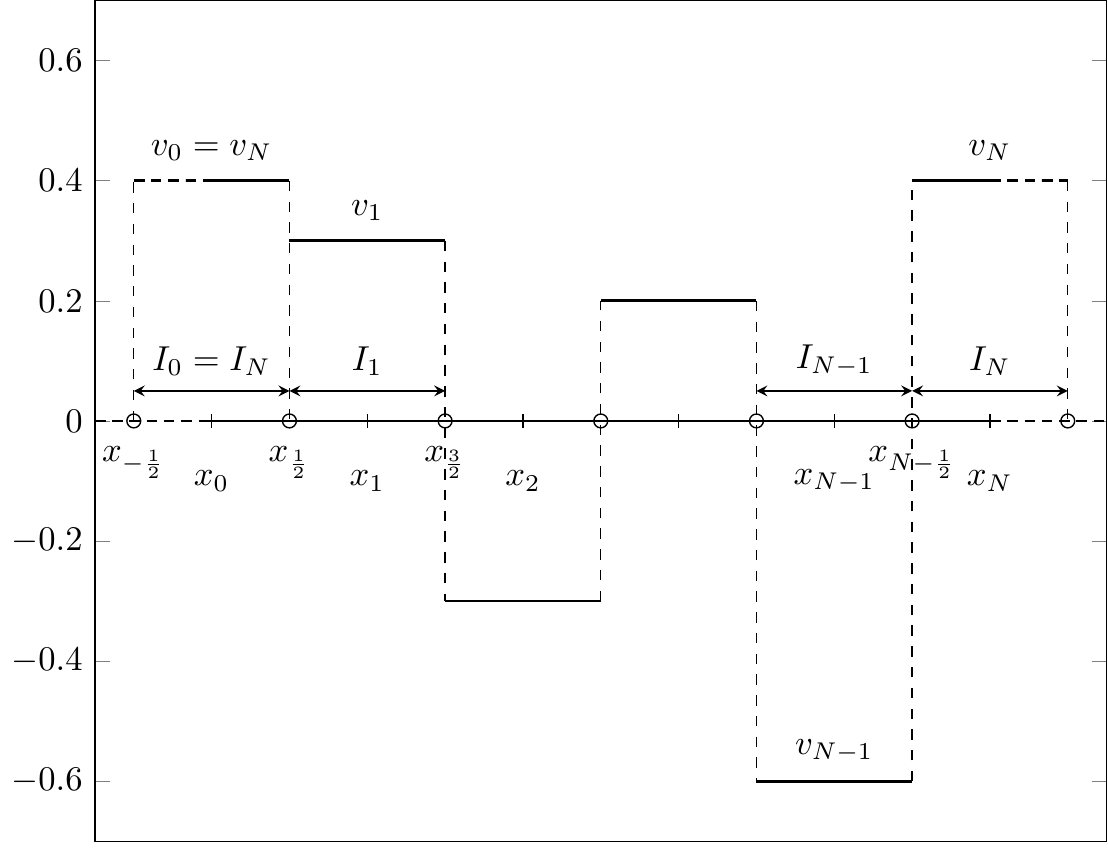}
\includegraphics[clip,scale=.7]{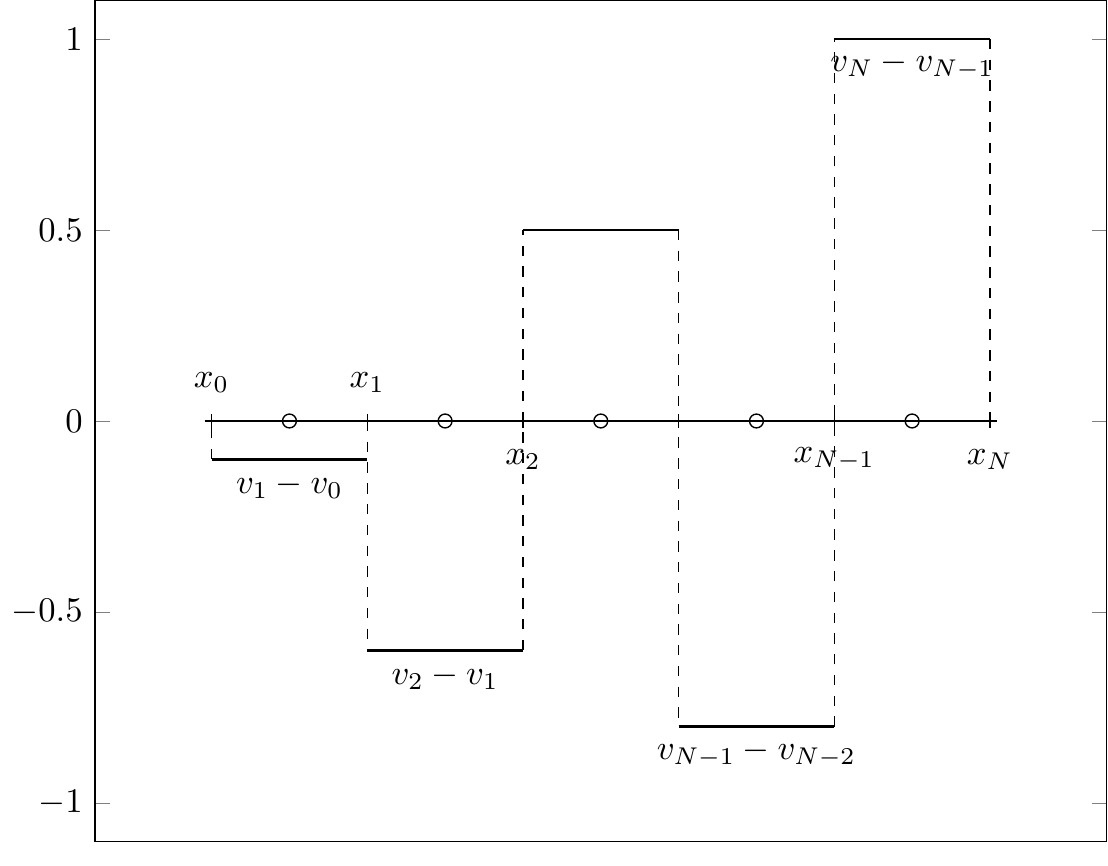}
\caption{An example of $v_h\in V_{h0}$ and $D_hv_h\in \widehat{V}_h$.}
\label{Fig:vh_and_Dhvh}
\end{figure}

Here we introduce the discretized problem for (P$1$);
\begin{equation}
(\mathrm{P}1)_h\quad \displaystyle\mathop{\mathrm{minimize}}_{u_h\in V_{h0}}\left\{\|\textbf{d}\|_1 + \dfrac{\lambda}{2}\|u_h-f_h\|_{H^{-1}_{\mathrm{av}}(\mathbb{T})}^2 : d_h = D_h u_h\in \widehat{V}_h\right\},
\end{equation} 
where $f_h\in V_{h0}$ is given data or $f_h=u_h^k$, and $\textbf{d} = (d_1,\dots, d_N)^{\mathrm{T}}$ for $d_h=\sum_{n=1}^Nd_n\boldsymbol{1}_{[x_{n-1},x_n)}$. 
Furthermore, we introduce the unconstrained problem
\begin{equation}
(\mathrm{P}2)_h\quad \displaystyle\mathop{\mathrm{minimize}}_{u_h\in V_{h0}, d_h\in\widehat{V}_h}\left\{\|\textbf{d}\|_1 + \dfrac{\lambda}{2}\|u_h-f_h\|_{H^{-1}_{\mathrm{av}}(\mathbb{T})}^2 + \dfrac{\mu}{2}\|d_h - D_h u_h\|_{L^2(I)}^2\right\}.
\end{equation}
\begin{remark}
In this paper, we use $\|d_h-D_hu_h\|_{L^2(I)}^2$. This enables to apply the shrinking method to minimization problem in the split Bregman framework.
\end{remark}

\subsection{Corresponding matrix form}
We reduce $(\mathrm{P}2)_h$ to the matrix formulation. Let $\textbf{d} = (d_1,\dots, d_N)^{\mathrm{T}}\in\mathbb{R}^N$ for $d_h=\sum_{n=1}^Nd_n\boldsymbol{1}_{[x_{n-1},x_n)}$, $\widetilde{\textbf{u}} = (u_1,\dots, u_N)^{\mathrm{T}}\in\mathbb{R}^N$ and $\textbf{u}=(u_1,\dots,u_{N-1})^{\mathrm{T}}\in\mathbb{R}^{N-1}$ for $u_h = \sum_{n=1}^Nu_n\boldsymbol{1}_{I_n}\in V_{h0}$, then we have
\begin{equation}
d_h-D_hu_h = \displaystyle\sum_{n=1}^N(d_n-(u_n-u_{n-1}))\boldsymbol{1}_{[x_{n-1},x_n)} = (\textbf{d}-S\widetilde{\textbf{u}})\cdot(\boldsymbol{1}_{[x_0,x_1)},\dots, \boldsymbol{1}_{[x_{N-1},x_N)})^{\mathrm{T}},
\end{equation}
where $S_N=h\nabla_h\in\mathbb{R}^{N\times N}$.
Furthermore, $u_h\in V_{h0}$ implies $u_N = -\sum_{n=1}^{N-1}u_n$, that is, $\widetilde{\textbf{u}} = R_N\textbf{u}$, where
\begin{equation}
\label{Def:R}
R_N = \begin{pmatrix}
1&0&\dots&0&\\
0&1&\dots&0&\\
\vdots&&\ddots&\vdots\\
0&0&\dots&1\\
-1&-1&\dots&-1
\end{pmatrix}\in\mathbb{R}^{N\times(N-1)}.
\end{equation}
Therefore
\begin{equation}
\label{Eq:L2_part}
\dfrac{\mu}{2}\|d_h-D_hu_h\|_{L^2(I)}^2 
= \dfrac{\mu h}{2}\|\textbf{d}-S_NR_N\textbf{u}\|_2^2.
\end{equation}

Next, we consider two expressions of $\|v_h\|_{H^{-1}_{\mathrm{av}}(\mathbb{T})}^2$ for $v_h\in V_{h0}$. Recall that equation \eqref{Lem:H-1av_norm} implies
\begin{equation}
\|v_h\|_{H^{-1}_{\mathrm{av}}(\mathbb{T})}^2 = \|\nabla(-\Delta_{\mathrm{av}})^{-1}v_h\|_{L^2(\mathbb{T})}^2.
\end{equation}
We propose two schemes for considering $\nabla(-\Delta_{\mathrm{av}})^{-1}$. The first scheme is to approximate $\nabla(-\Delta_{\mathrm{av}})^{-1}$ by using the discrete gradient $\nabla_h\in\mathbb{R}^{N\times N}$ and the discrete Laplacian
\begin{equation}
-\Delta_h = \nabla_h^{\mathrm{T}}\nabla_h = h^{-2}S_N^{\mathrm{T}}S_N = h^{-2}
\begin{pmatrix}
2&-1&0&\dots&-1\\
-1&2&-1&\dots&0\\
\vdots&&\ddots&&\vdots\\
-1&0&0&\dots&2
\end{pmatrix}\in\mathbb{R}^{N\times N}.
\end{equation}
Let $\widetilde{\textbf{v}} = (v_1,v_2,\dots, v_{N})^{\mathrm{T}}\in\mathbb{R}^N$ and $\textbf{v}=(v_1,\dots,v_{N-1})^{\mathrm{T}}\in\mathbb{R}^{N-1}$ for $v_h\in V_{h0}$, then $\widetilde{\textbf{v}}=R_N\textbf{v}$. 
We define $\textbf{w}\in\mathbb{R}^{N-1}$ and $\widetilde{\textbf{w}}\in\mathbb{R}^N$ for $w_h\in V_{h0}$ in the same way. Letting $\widetilde{\textbf{v}}= -\Delta_h\widetilde{\textbf{w}}$ implies
\begin{equation}
R_N\textbf{v} = -\Delta_hR_N\textbf{w}.
\end{equation}
Multiplying the (unique) pseudo-inverse matrix
\begin{equation}
L_N = \dfrac{1}{N}\begin{pmatrix}
N-1&-1&\dots&-1&-1\\
-1&N-1&\dots&-1&-1\\
\vdots&&\ddots&&\vdots\\
-1&-1&\dots&N-1&-1
\end{pmatrix}\in\mathbb{R}^{(N-1)\times N}
\end{equation}
yields $\textbf{v} = L_N(-\Delta_h)R_N\textbf{w} = h^{-2}L_NS_N^{\mathrm{T}}S_NR_N\textbf{w}$. For simplicity of notation, we let
\begin{subequations}
\begin{align}
A_N &= L_NS_N^{\mathrm{T}}S_NR_N,\\
(-\Delta_{\mathrm{av}})_h &= h^{-2}A_N = L_N(-\Delta_h)R_N.
\end{align}
\end{subequations}
It is easy to check that
\begin{equation}
A_N = \begin{pmatrix}
3&0&1&1&\dots&1&1&1\\
-1&2&-1&0&\dots&0&0&0\\
0&-1&2&-1&\dots&0&0&0\\
\vdots&&&\ddots&&&\vdots\\
0&0&0&0&\dots&-1&2&-1\\
1&1&1&1&\dots&1&0&3
\end{pmatrix}
\in \mathbb{R}^{(N-1)\times(N-1)}.
\end{equation}
satisfies $\det A_N=N^2\neq0$, therefore we have $\det(-\Delta_{\mathrm{av}})_h^{-1} \neq0$. 
This implies 
\begin{equation}
\left\{\begin{array}{rl}
\widetilde{\textbf{v}} &= -\Delta_h\widetilde{\textbf{w}},\\
\widetilde{\textbf{w}}&= R_N(-\Delta_{\mathrm{av}})_h^{-1}\textbf{v}.
\end{array}\right.
\end{equation}
Our first scheme is to approximate $(-\Delta_{\mathrm{av}})^{-1}$ by $R_N(-\Delta_{\mathrm{av}})_h^{-1}$, instead of $(-\Delta_h)^{-1}$ which does not exist. This yields
\begin{equation}
\nabla(-\Delta_{\mathrm{av}})^{-1}v_h\approx(\nabla_hR_N(-\Delta_{\mathrm{av}})_h^{-1}\textbf{v})\cdot(\boldsymbol{1}_{[x_0,x_1)},\dots,\boldsymbol{1}_{[x_{N-1},x_N)})^{\mathrm{T}},
\end{equation}
that is,
\begin{align*}
\|\nabla(-\Delta_{\mathrm{av}})^{-1}v_h\|_{L^2(I)}^2 &\approx \left\|\left(\nabla_hR_N(-\Delta_{\mathrm{av}})_h^{-1}\textbf{v}\right)\cdot(\boldsymbol{1}_{[x_0,x_1)},\dots,\boldsymbol{1}_{[x_{N-1},x_N)})\right\|_{L^2(I)}^2\\
&= h\|\nabla_hR_N(-\Delta_{\mathrm{av}})_h^{-1}\textbf{v}\|_2^2\\
&= h^3\|S_NR_NA_N^{-1}\textbf{v}\|_2^2.
\end{align*}
For simplicity, we let $J = SR_NA_N^{-1}\in\mathbb{R}^{N\times(N-1)}$, then our first scheme can be described as
\begin{equation}
\label{Eq:H-1_part1}
\dfrac{\lambda}{2}\|v_h\|_{H^{-1}_{\mathrm{av}}(\mathbb{T})}^2 \approx \dfrac{\lambda h^3}{2}\|J\textbf{v}\|_2^2
\end{equation}
for all $v_h\in V_{h0}$.
\begin{remark}
When we apply $H^{-s}_{\mathrm{av}}(\mathbb{T})$ norm $(0<s<1)$ to $(\mbox{P}1)_h$, the discrete inverse Laplacian $(-\Delta_{\mathrm{av}})_h^{-s}$ can be introduced by the discrete Fourier transform (for example, see \cite{GMR19}).
\end{remark}

\begin{figure}[tb]
\centering
\includegraphics[scale=.8,clip]{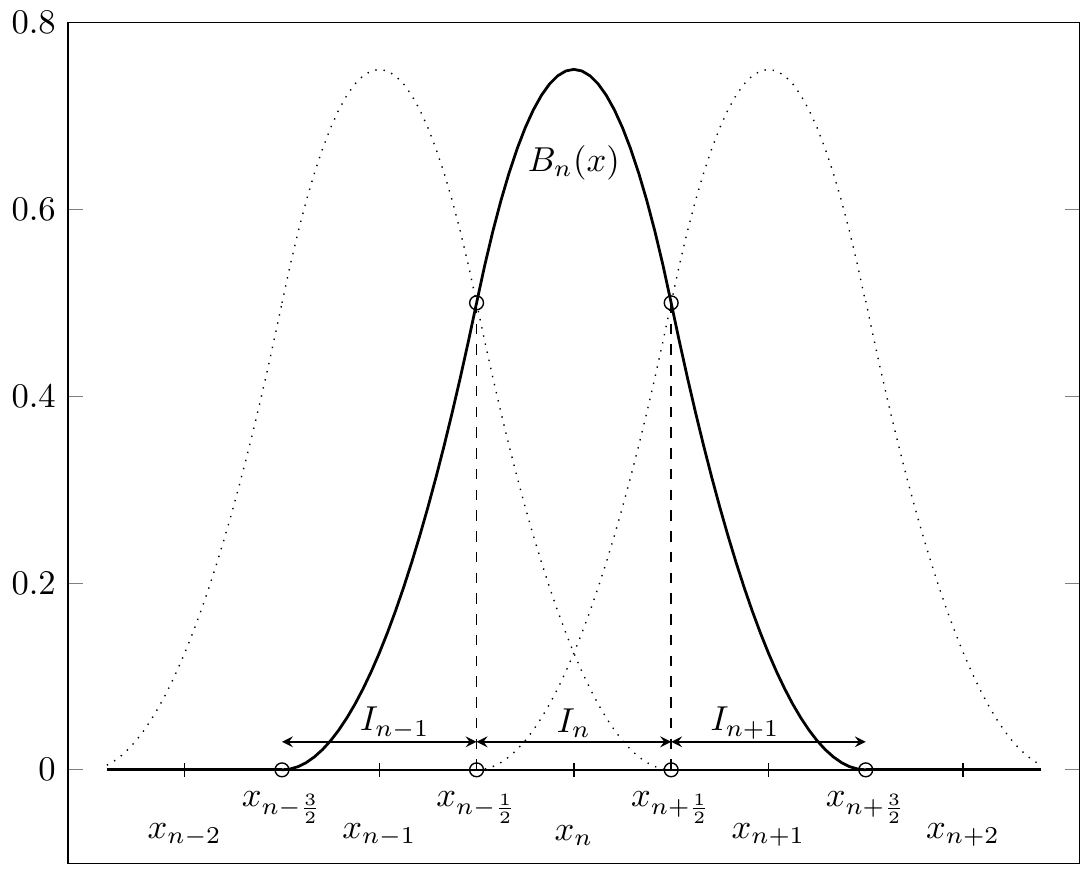}
\caption{The second degree B-spline basis functions}
\label{Fig:Bspline}
\end{figure}

Our second scheme is to compute $\nabla(-\Delta_{\mathrm{av}})^{-1}$ directly. It requires the second degree piecewise polynomial which has continuous derivative. We define the second degree periodic B-spline basis functions (see Figure \ref{Fig:Bspline})
\begin{equation}
B_n(x) = \left\{\begin{array}{rl}
\dfrac{(x-x_{n-\frac{3}{2}})^2}{2h^2}&\mbox{if }x\in I_{n-1},\\
\dfrac{(x-x_{n-\frac{1}{2}})(x_{n+\frac{1}{2}}-x)}{h^2}+\dfrac{1}{2}&\mbox{if }x\in I_n,\\
\dfrac{(x_{n+\frac{3}{2}}-x)^2}{2h^2}&\mbox{if }x\in I_{n+1},\\
0&\mbox{otherwise}.
\end{array}\right.
\end{equation}
We identify $B_{-1}\equiv B_{N-1}$, $B_0\equiv B_N$ and $B_1\equiv B_{N+1}$. The B-spline basis functions have continuous derivative (see Figure \ref{Fig:Deriv_Bspline})
\begin{equation}
\nabla B_n(x) =\left\{\begin{array}{rl}
(x-x_{n-\frac{3}{2}})h^{-2}&\mbox{if }x\in I_{n-1},\\
2(x_n-x)h^{-2}&\mbox{if }x\in I_n,\\
-(x_{n+\frac{3}{2}}-x)h^{-2}&\mbox{if }x\in I_{n+1},\\
0&\mbox{otherwise}.
\end{array}\right.
\end{equation}
Therefore we have
\begin{equation}
-\Delta B_n(x) = \left\{\begin{array}{rl}
-h^{-2}&\mbox{if }x\in I_{n-1},\\
2h^{-2}&\mbox{if }x\in I_n,\\
-h^{-2}&\mbox{if }x\in I_{n+1},\\
0&\mbox{otherwise}.
\end{array}\right.
\end{equation}
Fix $v_h\in V_{h0}$ arbitrarily, then there exits $w_h\in \operatorname{span}\{B_1,\dots, B_N\}$ such that $w_h = (-\Delta_{\mathrm{av}}^{-1})v_h\in H^1_{\mathrm{av}}(\mathbb{T})$. It is easy to check that 
\begin{equation}
\label{Eq:L1_Bspline}
\displaystyle\int_{\mathbb{T}}B_n(x)~dx = h\mbox{ for all }n=1,2,\dots, N.
\end{equation}
Let $\sum_{n=1}^Nw_n=0$, then equation \eqref{Eq:L1_Bspline} implies $\sum_{n=1}^Nw_nB_n\in H^1_{\mathrm{av}}(\mathbb{T})$.
Furthermore, we let
\begin{equation}
w_h = \displaystyle\sum_{n=1}^N w_nB_n\in H^1_{\mathrm{av}}(\mathbb{T}),\ \widetilde{\textbf{w}}=(w_1,\dots,w_N)^{\mathrm{T}}\in\mathbb{R}^N \mbox{ and } \textbf{w}=R_N\widetilde{\textbf{w}}\in\mathbb{R}^{N-1}.
\end{equation}
Then we have
\begin{equation}
v_h=-\Delta w_h = \displaystyle\sum_{n=1}^Nw_n (-\Delta B_n) = h^{-2}\sum_{n=1}^N(-w_{n-1}+2w_n-w_{n+1})\boldsymbol{1}_{I_n}\in V_{h0}.
\end{equation}
This implies
\begin{equation}
R_N\textbf{v}=\widetilde{\textbf{v}}= -\Delta_h\widetilde{\textbf{w}}=-\Delta_hR_N\textbf{w}=h^{-2}S_N^{\mathrm{T}}S_NR_N\textbf{w}.
\end{equation}
Multiplying the pseudo-inverse matrix $L_N$ yields
\begin{equation}
\textbf{v} =(-\Delta_{\mathrm{av}})_h\textbf{w}= h^{-2}A_N\textbf{w}.
\end{equation}
Therefore we have
\begin{equation}
\textbf{w} = (-\Delta_{\mathrm{av}})_h^{-1}\textbf{v}=h^2A_N^{-1}\textbf{v}.
\end{equation}
\begin{figure}[tb]
\centering
\includegraphics[scale=.7,clip]{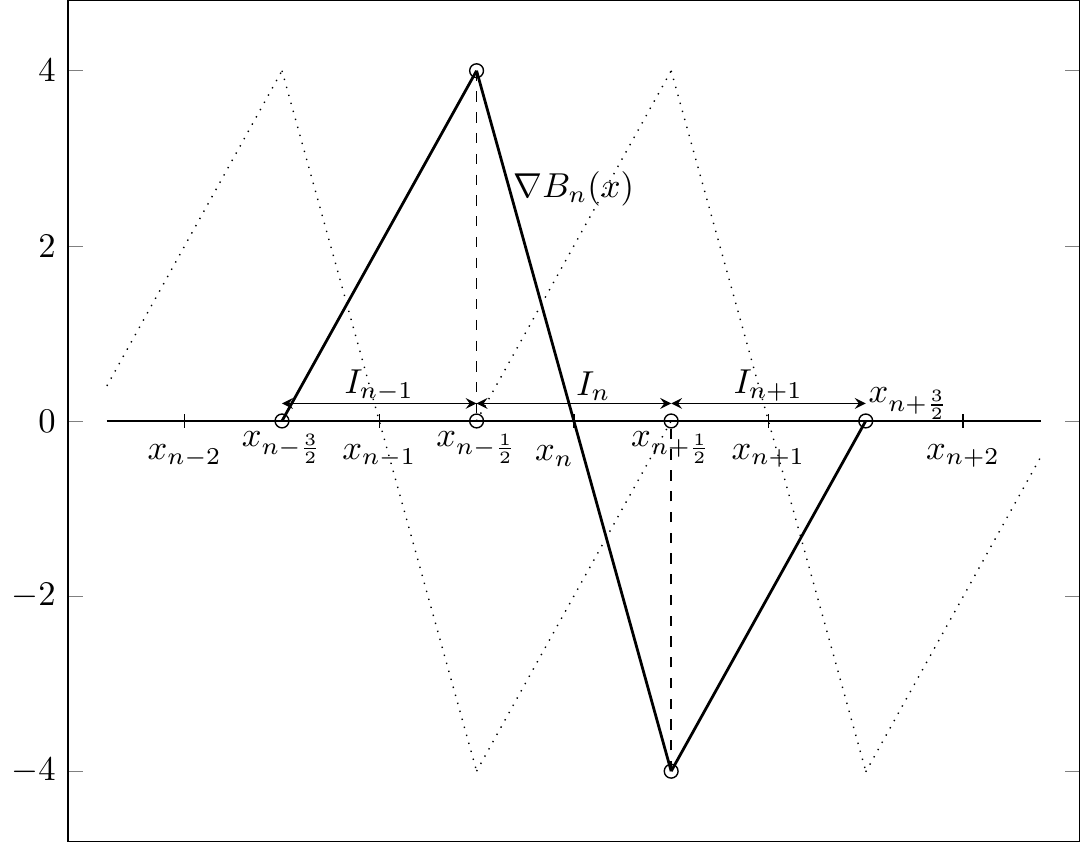}
\includegraphics[scale=.7,clip]{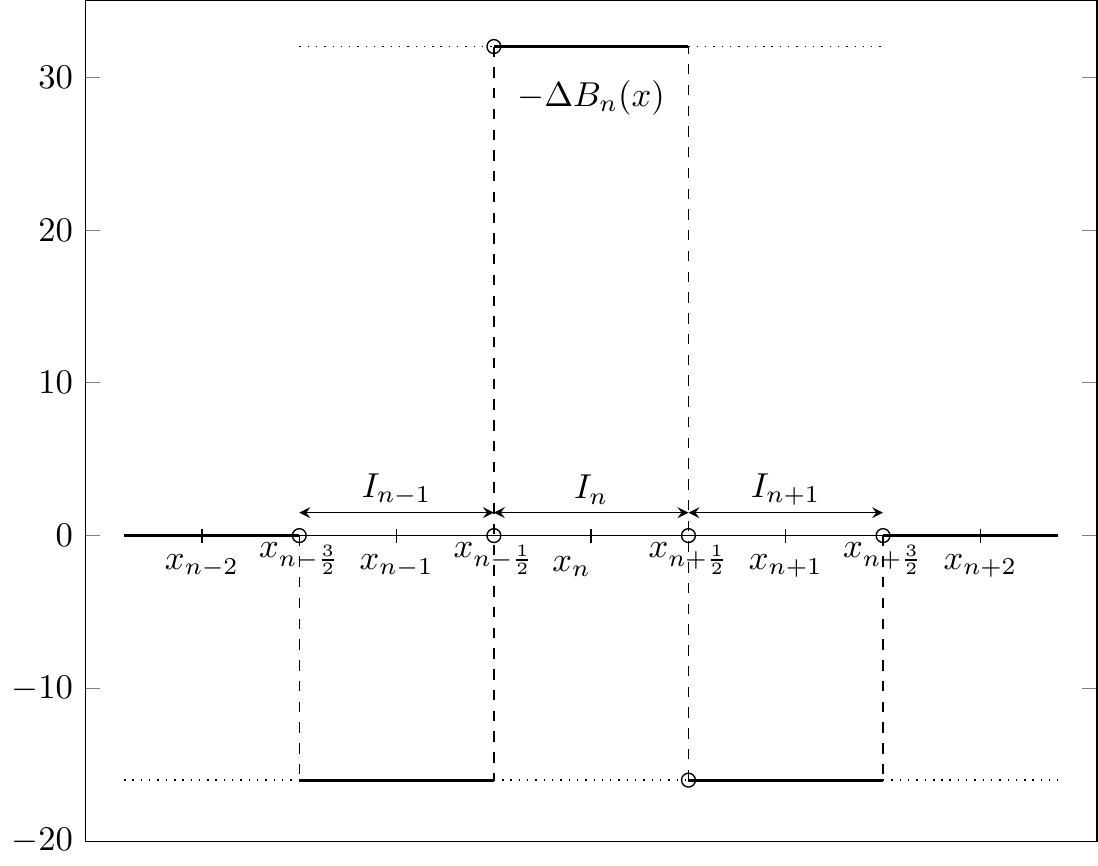}
\caption{The derivative of second degree B-spline basis functions.}
\label{Fig:Deriv_Bspline}
\end{figure}

The definition, combining with equation \eqref{Lem:H-1av_norm} gives $\|v_h\|_{H^{-1}_{\mathrm{av}}(\mathbb{T})}^2 = \|\nabla w_h\|_{L^2(\mathbb{T})}^2$, where
\begin{equation}
\nabla w_h = \displaystyle\sum_{n=1}^Nw_n\nabla B_n
\end{equation}
is a piecewise linear function which satisfies
\begin{equation}
\nabla w_h(x_{n-1/2}) = (w_n-w_{n-1})h^{-1}\mbox{ for all }n=1,\dots, N.
\end{equation}
This implies
\begin{equation}
\nabla w_h = \displaystyle\sum_{n=1}^N(w_n-w_{n-1})h^{-1}\phi_{n-1/2} = (\nabla_h\widetilde{\textbf{w}})\cdot(\phi_{1/2},\dots,\phi_{N-1/2})^{\mathrm{T}},
\end{equation}
where
\begin{equation}
\phi_{n-1/2}(x) = \left\{\begin{array}{rl}
(x-x_{n-3/2})h^{-1}&\mbox{if }x\in I_{n-1},\\
(x_{n+1/2}-x)h^{-1}&\mbox{if }x\in I_n,\\
0&\mbox{otherwise}.
\end{array}\right.
\end{equation}
We identify $\phi_{-1/2}=\phi_{N-1/2}$, $\phi_{1/2}=\phi_{N+1/2}$ (see Figure \ref{Fig:tent}). It is easy to check that
\begin{equation}
\displaystyle\int_{\mathbb{T}}\phi_{n-1/2}(x)\phi_{m-1/2}(x)~dx = \left\{
\begin{array}{rl}
2h/3&\mbox{if }n=m,\\
h/6&\mbox{if }|n-m|=1,\\
0&\mbox{otherwise}
\end{array}\right.\label{Eq:mass_mat_phi}
\end{equation}
for all $n=1,\dots,N$. Therefore we have
\begin{align*}
\|v_h\|_{H^{-1}_{\mathrm{av}}(\mathbb{T})}^2 &= \|\nabla w_h\|_{L^2(\mathbb{T})}^2\\
&= (\nabla_h\widetilde{\textbf{w}})^{\mathrm{T}}\begin{pmatrix}
2h/3&h/6&0&\dots&0&h/6\\
h/6&2h/3&h/6&\dots&0&0\\
\vdots&&\ddots&&&\vdots\\
h/6&0&0&\dots&h/6&2h/3
\end{pmatrix}\nabla_h\widetilde{\textbf{w}}\\
&= \dfrac{1}{h}(S_NR_N\textbf{w})^{\mathrm{T}}M_NS_NR_N\textbf{w},
\end{align*}
where
\begin{equation}
M_N = \begin{pmatrix}
2/3&1/6&0&\dots&0&1/6\\
1/6&2/3&1/6&\dots&0&0\\
\vdots&&\ddots&&&\vdots\\
1/6&0&0&\dots&1/6&2/3
\end{pmatrix}\in\mathbb{R}^{N\times N}.
\end{equation}
Let 
\begin{equation}
T=\begin{pmatrix}
a&0&\dots&0&b\\
b&a&\dots&0&0\\
\vdots&&\ddots&&\\
0&0&\dots&b&a\\
\end{pmatrix}\in\mathbb{R}^{N\times N},
\end{equation}
where $a=\dfrac{\sqrt{3}+1}{2\sqrt{3}}$ and $b = \dfrac{\sqrt{3}-1}{2\sqrt{3}}$, then $T^{\mathrm{T}}T=M_N$. Summarizing the above argument, our second scheme can be described as
\begin{align*}
\dfrac{\lambda}{2}\|v_h\|_{H^{-1}_{\mathrm{av}}(\mathbb{T})}^2 &= \dfrac{\lambda}{2h}(S_NR_N\textbf{w})^{\mathrm{T}}M_NS_NR_N\textbf{w}\\
&= \dfrac{\lambda}{2h}(S_NR_N(-\Delta_{\mathrm{av}})_h^{-1}\textbf{v})^{\mathrm{T}}T^{\mathrm{T}}TS_NR_N(-\Delta_{\mathrm{av}})_h^{-1}\textbf{v}\\
&=\dfrac{\lambda h^3}{2}\|TS_NR_NA_N^{-1}\textbf{v}\|_2^2.
\end{align*}
Let $H = TSR_NA_N^{-1} = TJ\in\mathbb{R}^{N\times (N-1)}$ for simplicity of notation, then we have
\begin{equation}
\label{Eq:H-1_part2}
\dfrac{\lambda}{2}\|v_h\|_{H^{-1}_{\mathrm{av}}(\mathbb{T})}^2 = \dfrac{\lambda h^3}{2}\|H\textbf{v}\|_2^2.
\end{equation}

Applying equation \eqref{Eq:L2_part}, \eqref{Eq:H-1_part1} and \eqref{Eq:H-1_part2} to $(\mbox{P}2)_h$ implies the following two discretized problems;
\begin{subequations}
\begin{align}
&\displaystyle\mathop{\mathrm{minimize}}_{\textbf{u}\in\mathbb{R}^{N-1}, \textbf{d}\in\mathbb{R}^N}\left\{\|\textbf{d}\|_1+\dfrac{\lambda h^3}{2}\|J(\textbf{u}-\textbf{f})\|_2^2+\dfrac{\mu h}{2}\|\textbf{d}-S_NR_N\textbf{u}\|_2^2\right\},\label{Scheme:SB1}\\
&\displaystyle\mathop{\mathrm{minimize}}_{\textbf{u}\in\mathbb{R}^{N-1}, \textbf{d}\in\mathbb{R}^N}\left\{\|\textbf{d}\|_1+\dfrac{\lambda h^3}{2}\|H(\textbf{u}-\textbf{f})\|_2^2+\dfrac{\mu h}{2}\|\textbf{d}-S_NR_N\textbf{u}\|_2^2\right\}\label{Scheme:SB2},
\end{align}
\end{subequations}
where $\textbf{f}\in\mathbb{R}^{N-1}$ is given as $f\in V_{h0}$ or $\textbf{f} = \textbf{u}^k$. Recall that the matrix $J$ is introduced by the approximation $\nabla(-\Delta_{\mathrm{av}})^{-1}\approx \nabla_hR_N(-\Delta_{\mathrm{av}})_h^{-1}$. On the other hand, we obtain $H$ by using $\nabla(-\Delta_{\mathrm{av}})^{-1}$ exactly. Therefore \eqref{Scheme:SB1} can be regarded as an approximation of \eqref{Scheme:SB2}, which is equivalent to $(\mbox{P}2)_h$.

\begin{figure}[tb]
\centering
\includegraphics[scale=.8]{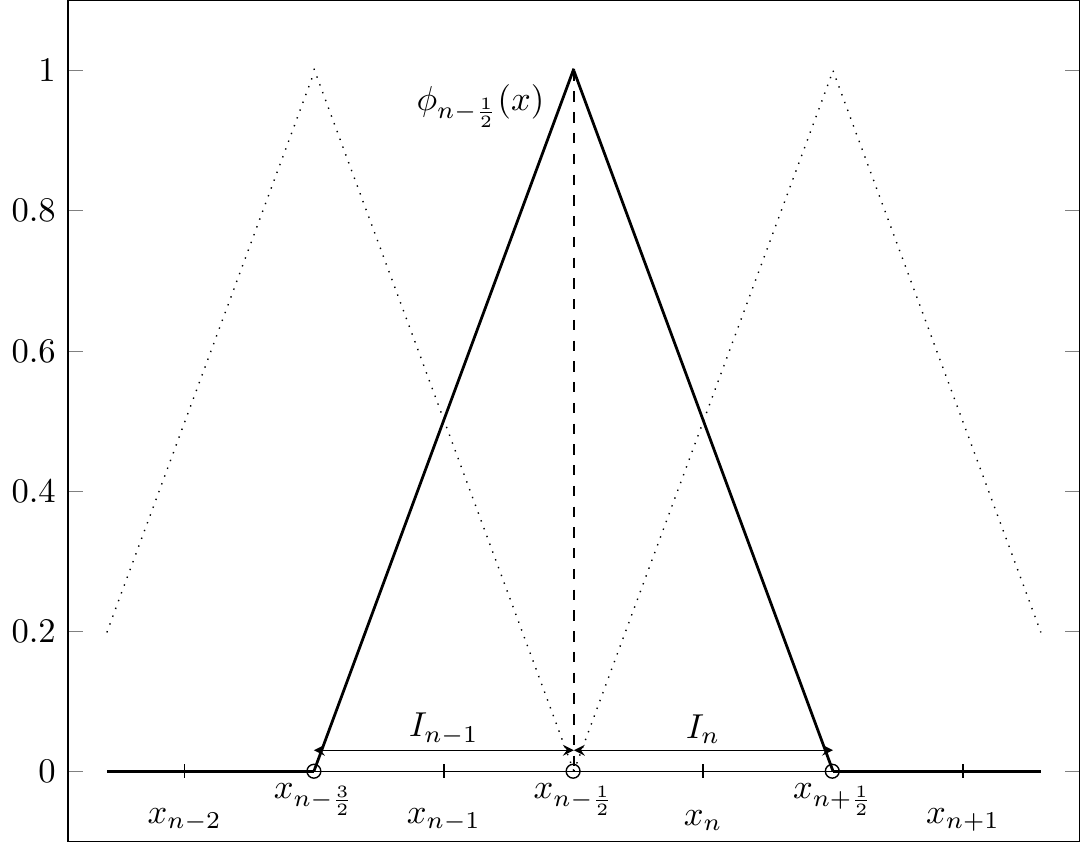}
\caption{The piecewise linear basis functions}
\label{Fig:tent}
\end{figure}

\section{Split Bregman framework}
\label{Sec:SplitBregman}
In this section, we review the alternating split Bregman framework in \cite{GO09} for the problem 
\begin{equation}
\label{Prob:Generalized}
(\mbox{P}3K)_h\quad \displaystyle\mathop{\mathrm{minimize}}_{\textbf{u}\in\mathbb{R}^{N-1}, \textbf{d}\in\mathbb{R}^N}\left\{\|\textbf{d}\|_1+\dfrac{\lambda h^3}{2}\|K(\textbf{u}-\textbf{f})\|_2^2+\dfrac{\mu h}{2}\|\textbf{d}-S_NR_N\textbf{u}\|_2^2\right\},
\end{equation}
where $K\in\mathbb{R}^{N\times(N-1)}$ is equal to $J$ or $H$. Recall that $(\mbox{P}3K)_h$ is an approximation of the discrete problem for $(\mbox{P}0)$;
\begin{equation}
(\mbox{P}0K)_h\quad \displaystyle\mathop{\mathrm{minimize}}_{\textbf{u}\in\mathbb{R}^{N-1}}\left\{\|S_NR_N\textbf{u}\|_1+\dfrac{\lambda h^3}{2}\|K(\textbf{u}-\textbf{f})\|_2^2\right\}.
\end{equation} 
Let
\begin{equation}
\Psi(\textbf{u},\textbf{d}) = \|\textbf{d}\|_1+\dfrac{\lambda h^3}{2}\|K(\textbf{u}-\textbf{f})\|_2.
\end{equation}
The Bregman method replaces $\Psi(\textbf{u},\textbf{d})$ into the Bregman distance and iteratively solves
\begin{equation}
\label{Eq:Bregman}
(\textbf{u}^{k+1},\textbf{d}^{k+1}) = \displaystyle\mathop{\mathrm{argmin}}_{\textbf{u}\in\mathbb{R}^{N-1}, \textbf{d}\in\mathbb{R}^N}\left\{D_{\Psi}^{\textbf{p}^k}((\textbf{u},\textbf{d}),(\textbf{u}^k,\textbf{d}^k)) + \dfrac{\mu h}{2}\|\textbf{d}-S_NR_N\textbf{u}\|_2^2\right\},
\end{equation}
where the Bregman distance $D_{\Psi}^{\textbf{p}^k}$ is defined as
\begin{equation}
D_{\Psi}^{\textbf{p}^k}((\textbf{u},\textbf{d}),(\textbf{u}^k,\textbf{d}^k)) = \Psi(\textbf{u},\textbf{d}) - \Psi(\textbf{u}^k,\textbf{d}^k)-\textbf{p}_u^k\cdot(\textbf{u}-\textbf{u}^k)-\textbf{p}_d^k\cdot(\textbf{d}-\textbf{d}^k),
\end{equation}
and $\textbf{p}^k = (\textbf{p}_u^k,\textbf{p}_d^k)\in \mathbb{R}^{N-1}\times\mathbb{R}^N$ is defined as
\begin{subequations}
\begin{align}
\textbf{p}_u^{k+1} &= \textbf{p}_u^k - \mu h(S_NR_N)^{\mathrm{T}}(S_NR_N\textbf{u}^{k+1}-\textbf{d}^{k+1})\mbox{ and }\textbf{p}_u^0 =\textbf{0}\in\mathbb{R}^{N-1}\\
\textbf{p}_d^{k+1} &= \textbf{p}_d^k - \mu h(\textbf{d}^{k+1}-S_NR_N\textbf{u}^{k+1})\mbox{ and }\textbf{p}_d^0=\textbf{0}\in\mathbb{R}^N.
\end{align}
\end{subequations}
Thanks to $\Psi:\mathbb{R}^{N-1}\times\mathbb{R}^N\to\mathbb{R}$ is convex and lower semi-continuous, the Bregman distance $D_{\Psi}^{\textbf{p}^k}(\cdot,(\textbf{u}^k,\textbf{d}^k))$ is also convex and lower semi-continuous. Applying the usual existence result of convex analysis (see \cite[Cor3.23]{Bre11}) gives that there exists a minimizer $(\textbf{u}^{k+1},\textbf{d}^{k+1})$. Furthermore, by using induction we can show that
\begin{equation}
\left(\Psi(\textbf{u},\textbf{d}) +\dfrac{\mu h}{2}\|\textbf{d}-S_NR_N\textbf{u}-\boldsymbol{\alpha}^{k}\|_2^2\right) - \left(D_{\Psi}^{\textbf{p}^k}((\textbf{u},\textbf{d}),(\textbf{u}^k,\textbf{d}^k)) + \dfrac{\mu h}{2}\|\textbf{d}-S_NR_N\textbf{u}\|_2^2\right)
\end{equation}
is independent of $(\textbf{u},\textbf{d})$, where $\boldsymbol{\alpha}^{k+1}\in\mathbb{R}^N$ is defined as
\begin{equation}
\boldsymbol{\alpha}^{k+1} = \boldsymbol{\alpha}^k -(\textbf{d}^{k+1}-S_NR_N\textbf{u}^{k+1})\mbox{ and }\boldsymbol{\alpha}^{0}=\textbf{0}.
\end{equation}
This implies the minimizer $(\textbf{u}^{k+1},\textbf{d}^{k+1})$ of problem \eqref{Eq:Bregman} satisfies
\begin{equation}
(\textbf{u}^{k+1},\textbf{d}^{k+1}) = \displaystyle\mathop{\mathrm{argmin}}_{\textbf{u}\in\mathbb{R}^{N-1}, \textbf{d}\in\mathbb{R}^N}\left\{\|\textbf{d}\|_1+\dfrac{\lambda h^3}{2}\|K(\textbf{u}-\textbf{f})\|_2^2+\dfrac{\mu h}{2}\|\textbf{d}-S_NR_N\textbf{u}-\boldsymbol{\alpha}^{k}\|_2^2\right\}.
\end{equation}
This is the split Bregman iteration for the problem $(\mbox{P}3K)_h$. 
Finally, we apply the alternating split Bregman algorithm and obtain
\begin{subequations}
\begin{empheq}[left=(\mbox{P}4K)_h\quad\empheqlbrace]{align}
\textbf{u}^{k+1} &= \displaystyle\mathop{\mathrm{argmin}}_{\textbf{u}\in\mathbb{R}^{N-1}}\left\{\dfrac{\lambda h^3}{2}\|K(\textbf{u}-\textbf{f})\|_2^2 + \dfrac{\mu h}{2}\|\textbf{d}^k-S_NR_N\textbf{u}-\boldsymbol{\alpha}^{k}\|_2^2\right\},\label{Eq:P4h2}\\
\textbf{d}^{k+1} &= \displaystyle\mathop{\mathrm{argmin}}_{\textbf{d}\in\mathbb{R}^N}\left\{\|\textbf{d}\|_1+\dfrac{\mu h}{2}\|\textbf{d}-S_NR_N\textbf{u}^{k+1}-\boldsymbol{\alpha}^{k}\|_2^2\right\}\label{Eq:P4h3},\\
\boldsymbol{\alpha}^{k+1} &= \boldsymbol{\alpha}^k-\textbf{d}^{k+1}+S_NR_N\textbf{u}^{k+1},
\end{empheq}
\end{subequations}
where $\textbf{f}\in\mathbb{R}^{N-1}$ is given data or $\textbf{f}=\textbf{u}^k$, $\boldsymbol{\alpha}^0 = \textbf{0}$, $\textbf{u}^0$ is given as $\textbf{0}$ or initial condition, and $\textbf{d}^0 = S_NR_N\textbf{u}^0$. This satisfies the following convergence result.
\begin{lemma}[Theorem 3.2 of \cite{COS10}]
Suppose that $(\mbox{P}0K)_h$ has a minimizer $\textbf{u}^*\in\mathbb{R}^{N-1}$, then $\textbf{u}^k$ which determined by $(\mbox{P}4K)_h$ satisfies
\begin{equation}
\lim_{k\to\infty}\|S_NR_N\textbf{u}^{k+1}\|_1 + \dfrac{\lambda h^3}{2}\|K(\textbf{u}^{k+1}-\textbf{f})\|_2^2 = \|S_NR_N\textbf{u}^*\|_1 + \dfrac{\lambda h^3}{2}\|K(\textbf{u}^*-\textbf{f})\|_2^2.
\end{equation}
Furthermore, if a minimizer $\textbf{u}^*$ of $(\mbox{P}0K)_h$ is unique, then $\lim_{k\to\infty}\|\textbf{u}^{k+1}-\textbf{u}^*\|_2=0$.
\end{lemma}

The functional in \eqref{Eq:P4h2} is differentiable with respect to $\textbf{u}$, and the minimization \eqref{Eq:P4h3} can be reduced to the shrinking method
\begin{equation}
(\textbf{d}^{k+1})_n = \operatorname{shrink}\left((S_NR_N\textbf{u}^{k+1}+\boldsymbol{\alpha}^{k})_n,\dfrac{1}{\mu h}\right),
\end{equation}
where $(\textbf{v})_n$ is the $n$-entry of vector $\textbf{v}$ and
\begin{equation}
\label{Op:shrink}
\operatorname{shrink}(\rho,a) = \dfrac{\rho}{|\rho|}\max\{|\rho|-a,0\}.
\end{equation}
Therefore, the problem $(\mbox{P}4K)_h$ introduces
\begin{equation*}
(\mbox{P}5K)_h\quad \left\{\begin{array}{rl}
\textbf{u}^{k+1} &= \left(\lambda h^3K^{\mathrm{T}}K+\mu h(S_NR_N)^{\mathrm{T}}S_NR_N\right)^{-1}\left(\lambda h^3K^{\mathrm{T}}K\textbf{f}+\mu h(S_NR_N)^{\mathrm{T}}(\textbf{d}^k-\boldsymbol{\alpha}^{k})\right),\\
(\textbf{d}^{k+1})_n &= \operatorname{shrink}\left((S_NR_N\textbf{u}^{k+1}+\boldsymbol{\alpha}^{k})_n,\dfrac{1}{\mu h}\right)\mbox{ for all }n=1,\dots, N,\\
\boldsymbol{\alpha}^{k+1} &= \boldsymbol{\alpha}^k-\textbf{d}^{k+1}+S_NR_N\textbf{u}^{k+1},
\end{array}\right.
\end{equation*}
where $\textbf{f}\in\mathbb{R}^{N-1}$ is given data or $\textbf{f}=\textbf{u}^k$, $\boldsymbol{\alpha}^0 = \textbf{0}$, $\textbf{u}^0$ is given as $\textbf{0}$ or initial condition, and $\textbf{d}^0 = S_NR_N\textbf{u}^0$.

\begin{figure}[tb]
\centering
\includegraphics[clip]{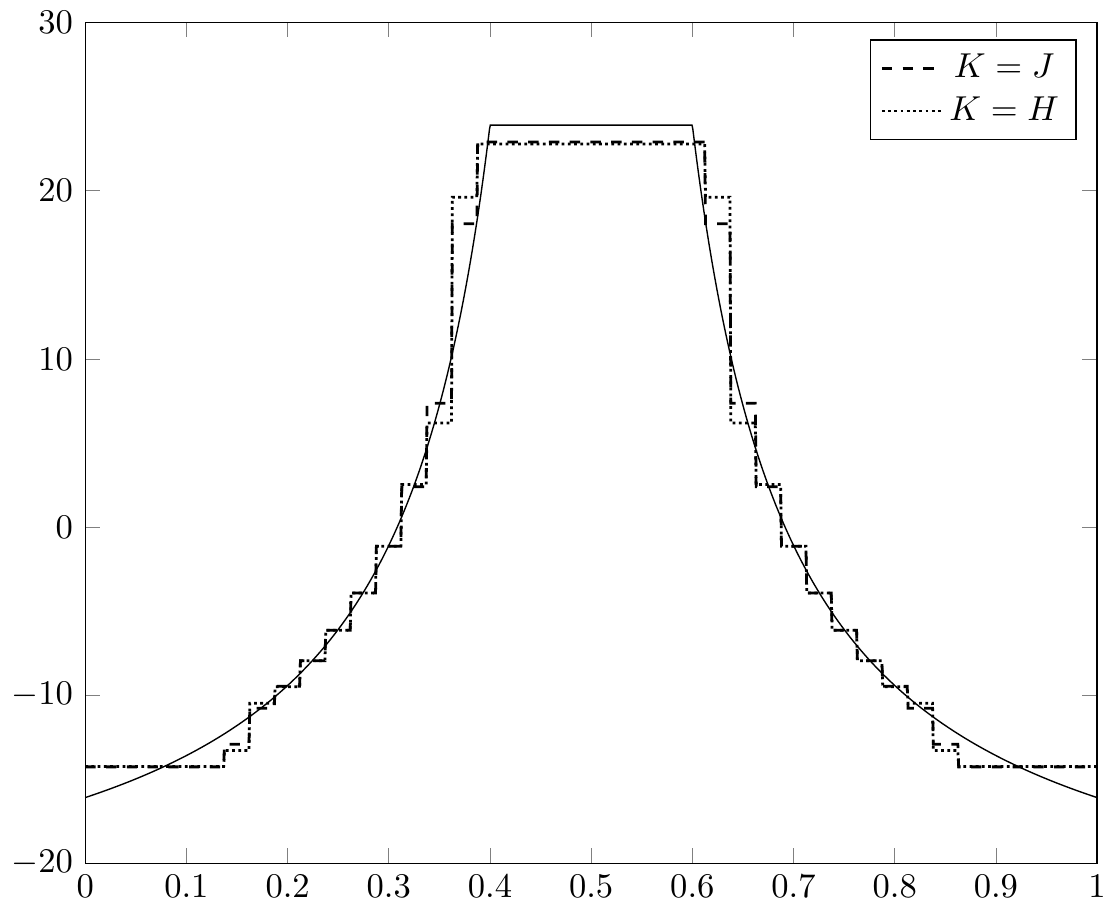}
\caption{The difference between $K=J$ and $K=H$.}
\label{Fig:example3}
\end{figure}

\section{Shrinking method for Spohn's model}
\label{Sec:Spohn}
In this section, we consider the split Bregman framework for Spohn's model
\begin{equation}
u_t = -\Delta\left(\operatorname{div}\left(\beta\dfrac{\nabla u}{|\nabla u|}+|\nabla u|^{p-2}\nabla u\right)\right),
\end{equation}
which can be regarded as the gradient flow problem for energy functional
\begin{equation}
\widetilde{\Phi}(u) = \beta\displaystyle\int_{\mathbb{T}}|Du| + \dfrac{1}{p}\int_{\mathbb{T}}|Du|^p,
\end{equation}
where $\beta>0$ and $p>1$. This energy is considered in model for the relaxation of a crystalline surface below the roughening temperature (for example, see \cite{KV10}). If we replace $w = (Du)^p$, the alternating split Bregman method introduces nonlinear problem. In this paper we always assume $p=3$, and we apply the constraint $d=Du$ to $|Du|^3$. The alternating split Bregman method implies
\begin{subequations}
\begin{empheq}[left=\empheqlbrace]{align}
\textbf{u}^{k+1} &= \displaystyle\mathop{\mathrm{argmin}}_{\textbf{u}\in\mathbb{R}^{N-1}}\left\{\dfrac{\lambda h^3}{2}\|K(\textbf{u}-\textbf{f})\|_2^2 + \dfrac{\mu h}{2}\|\textbf{d}^k-S_NR_N\textbf{u}-\boldsymbol{\alpha}^{k}\|_2^2\right\},\\
\textbf{d}^{k+1} &= \displaystyle\mathop{\mathrm{argmin}}_{\textbf{d}\in\mathbb{R}^N}\left\{\beta\|\textbf{d}\|_1+ \dfrac{1}{p}\|\textbf{d}\|_p^p +\dfrac{\mu h}{2}\|\textbf{d}-S_NR_N\textbf{u}^{k+1}-\boldsymbol{\alpha}^{k}\|_2^2\right\}\label{Eq:Crystalline_energy},\\
\boldsymbol{\alpha}^{k+1} &= \boldsymbol{\alpha}^k-\textbf{d}^{k+1}+S_NR_N\textbf{u}^{k+1},
\end{empheq}
\end{subequations}
We consider the Euler-Lagrange equation for equation \eqref{Eq:Crystalline_energy};
\begin{equation}
\beta\left(\dfrac{(\textbf{d}^{k+1})_n}{|(\textbf{d}^{k+1})_n|}\right)_{1\le n\le N} + \left((\textbf{d}^{k+1})_n|(\textbf{d}^{k+1})_n|^{p-2}\right)_{1\le n\le N} +\mu h(\textbf{d}^{k+1}-SR_N\textbf{u}^{k+1}-\boldsymbol{\alpha}^k)=0.
\end{equation}
For simplicity of notation, we let $x=(\textbf{d}^{k+1})_n$, $a=1/(\mu h)>0$ and $\rho = (SR_N\textbf{u}^{k+1}+\boldsymbol{\alpha}^k)_n$. This, combining with $p=3$ gives
\begin{equation}
\beta\dfrac{x}{|x|}+x|x|+\dfrac{1}{a}(x-\rho)=0.
\end{equation}
Suppose that $x>0$, then we have $a\beta<\rho$ and 
\begin{equation}
x=\dfrac{1}{2}\left(-\dfrac{1}{a}+\sqrt{\dfrac{1}{a^2}-4\left(\beta-\dfrac{\rho}{a}\right)}\right).
\end{equation}
By the similar way, supposing $x<0$ yields $\rho<-a\beta$ and
\begin{equation}
x = \dfrac{1}{2}\left(\dfrac{1}{a}-\sqrt{\dfrac{1}{a^2}-4\left(\beta+\dfrac{\rho}{a}\right)}\right).
\end{equation}
If $-a\beta<\rho<a\beta$, we let $x=0$. These observations provide the shrinking operator of the form
\begin{equation}
\label{Eq:Shrink_for_Spohn}
x = \dfrac{\rho}{2a|\rho|}\left(-1+\sqrt{1+4a\max\{|\rho|-a\beta,0\}}\right).
\end{equation}
Applying this to equation \eqref{Eq:Crystalline_energy} gives
\begin{equation*}
\left\{\begin{array}{rl}
\textbf{u}^{k+1} &= \left(\lambda h^3K^{\mathrm{T}}K+\mu h(S_NR_N)^{\mathrm{T}}S_NR_N\right)^{-1}\left(\lambda h^3K^{\mathrm{T}}K\textbf{f}+\mu h(S_NR_N)^{\mathrm{T}}(\textbf{d}^k-\boldsymbol{\alpha}^{k})\right),\\
(\textbf{d}^{k+1})_n &=\dfrac{\mu h\rho_n^{k+1}}{2|\rho_n^{k+1}|}\left(-1+\sqrt{1+\dfrac{4}{\mu h}\max\left\{|\rho_n^{k+1}|-\dfrac{\beta}{\mu h},0\right\}}\right) \mbox{ for all }n=1,\dots, N,\\
\boldsymbol{\alpha}^{k+1} &= \boldsymbol{\alpha}^k-\textbf{d}^{k+1}+S_NR_N\textbf{u}^{k+1},
\end{array}\right.
\end{equation*}
where $\rho_n^{k+1}=(S_NR_N\textbf{u}^{k+1}+\boldsymbol{\alpha}^k)_n$.

\begin{figure}[tb]
\centering
\subfloat[][$K=J$]{\includegraphics[clip,scale=.7]{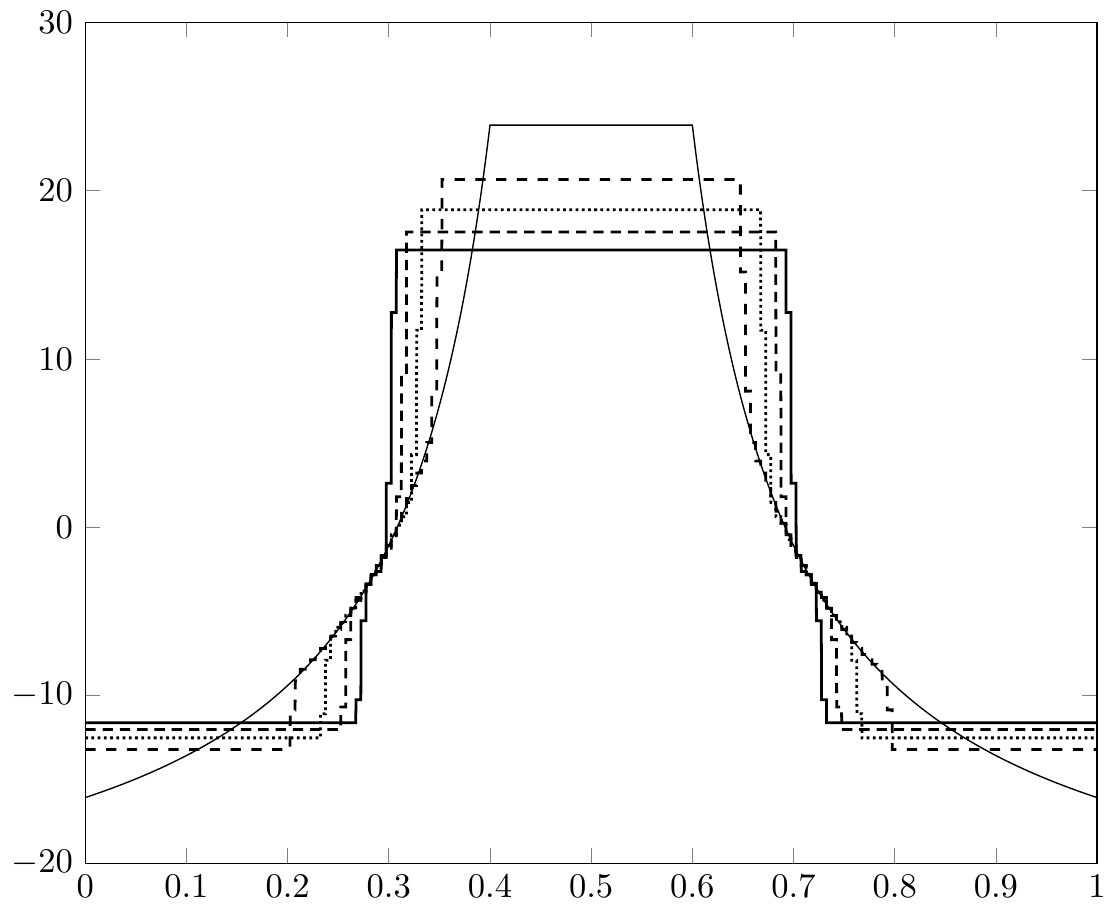}\label{Fig:K=J}}
\subfloat[][$K=H$]{\includegraphics[clip,scale=.7]{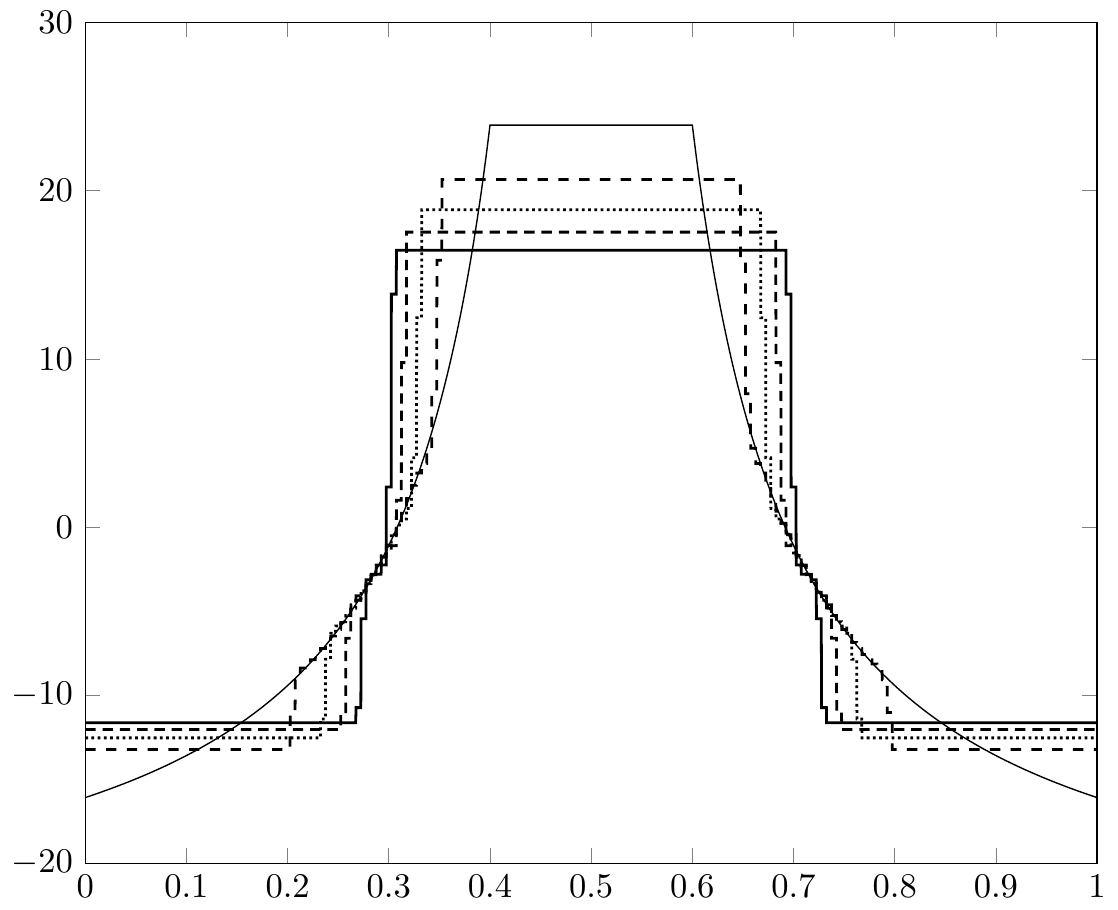}\label{Fig:K=H}}
\caption{Numerical examples of the gradient flow.}
\label{Fig:gradient_flow}
\end{figure}

\section{Numerical example}
\label{Sec:NumExample}
\subsection{Example 1: Comparison of two schemes}
Here we show numerical examples of $(\mbox{P}5K)_h$. Note that equation \eqref{Eq:P4h3} implies that $\mu$ should satisfy $\mu=O(h^{-1})$. Moreover, this and equation \eqref{Eq:P4h2} yield $\lambda = O(h^{-3})$ is necessary for reasonable computation.

In this paper, we always regard $\mathbb{T}$ as an interval $[0,1]$ with periodic boundary condition.
Our first numerical example is the gradient flow \eqref{Eq:GradientFlow} with the initial condition
\begin{equation}
u^0(x) = \left\{\begin{array}{ll}
10(4-\log5)&\mbox{if }|x-1/2|\le 1/10,\\
\dfrac{5}{|x-1/2|}-10(1+\log5)&\mbox{otherwise.}
\end{array}\right.
\end{equation}
Note that the similar example is computed in \cite{GMR19}. They essentially apply the matrix $J$ and compute the gradient flow problem without split Bregman method. Their scheme requires $\tau = \lambda^{-1} = O(h^5)$ for $H^{-1}_{\mathrm{av}}$ fidelity.

We check the difference between $K=J$ and $K=H$. Figure \ref{Fig:example3} shows two numerical results with the same parameters $N=40$, $\lambda=h^{-3}$ and $\mu=5h^{-1}$. 
Numerical results $\textbf{u}^k\in\mathbb{R}^{N-1}$ are represented as piecewise constant functions $u_h^k\in V_{h0}$.
They are different because the matrix $J$ is introduced by discrete gradient and discrete inverse Laplacian. This difference is expected to be small if we consider sufficiently small $h$. 
Figure \ref{Fig:gradient_flow} shows evolution of numerical solutions for $N=200$, $\lambda=h^{-3}$ and $\mu=5h^{-1}$. We infer from them that \eqref{Scheme:SB1} can provide sufficiently accurate result.

\subsection{Example 2: Discontinuity and symmetry}
Our second numerical example for \eqref{Eq:GradientFlow} is
\begin{equation}
\label{Eq:initial_data}
u^0(x) = \left\{\begin{array}{ll}
-a(1/4-r)^3&\mbox{if }0<x<r\mbox{ or }1-r<x<1,\\
a(x-1/4)^3&\mbox{if }r<x<1/2-r,\\
a(1/4-r)^3&\mbox{if }1/2-r<x<1/2+r\\
-a(x-3/4)^3&\mbox{if }1/2+r<x<1-r,
\end{array}\right.
\end{equation}
where $a = 450$ and $r=1/15$. 
In \cite{GG10}, a class of initial data including \eqref{Eq:initial_data} as an example has been studied analytically. They rigorously proved that the solution becomes discontinuous instantaneously. Their analysis gives an exact profile of the fourth order gradient flow. Note that because of uniqueness of a solution, the symmetry of initial profile is preserved during evolution.
We can check that our numerical result shows the discontinuity and symmetry approximately (see Figure \ref{Fig:grad41}, \ref{Fig:grad42}). We use $K=J$, $N=200$, $\lambda=25h^{-3}$ and $\mu=15h^{-1}$. Furthermore, we note that we can compute until $\textbf{u}^k\approx \textbf{0}$ easily, because our scheme can be stable for $\tau=\lambda^{-1}=O(h^3)$. 

\begin{figure}[tb]
\centering
\subfloat[][Evolution of numerical result]{\includegraphics[clip,scale=.55]{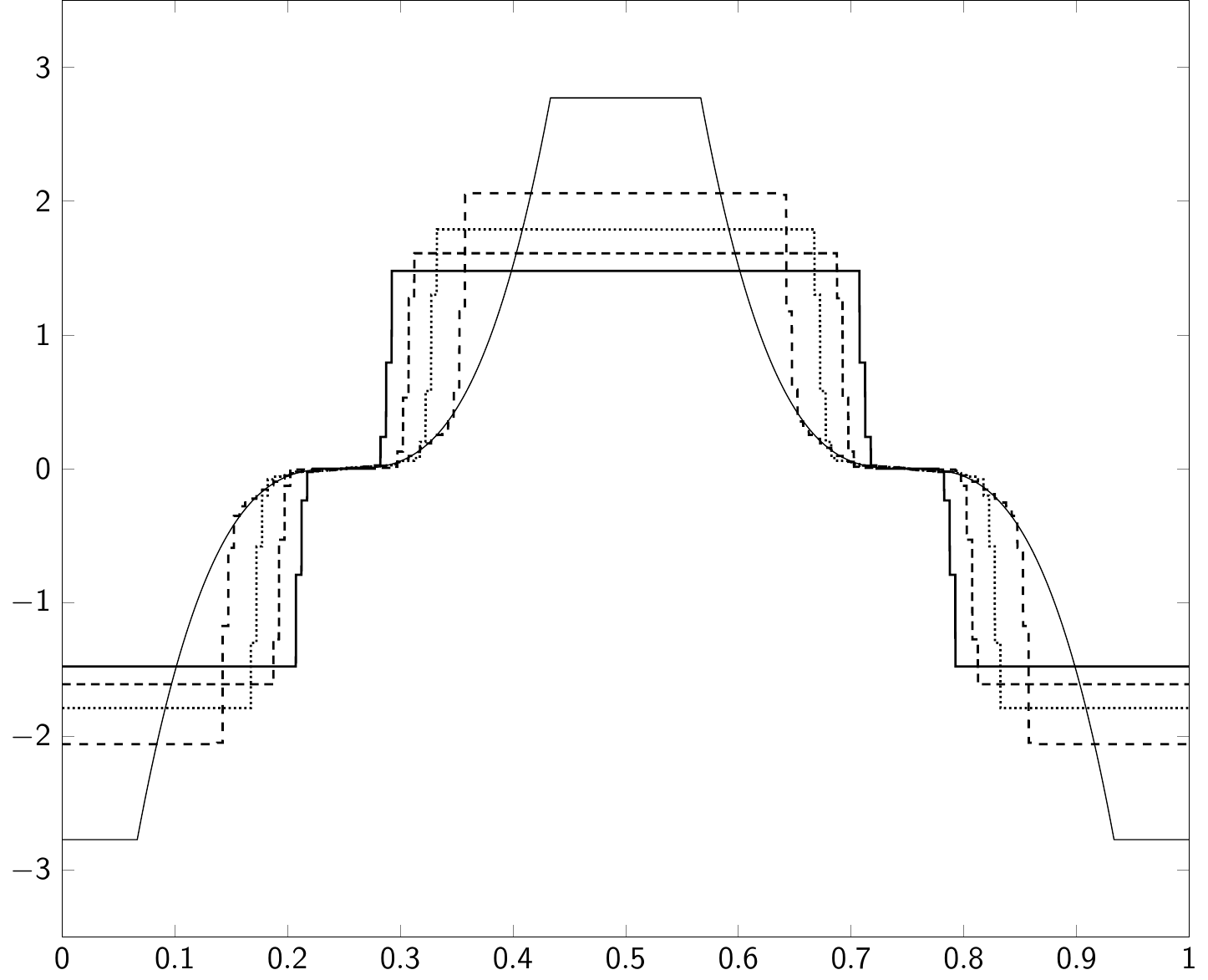}\label{Fig:grad41}}
\subfloat[][Evolution of numerical result around $u^k\approx0$]{\includegraphics[clip,scale=.55]{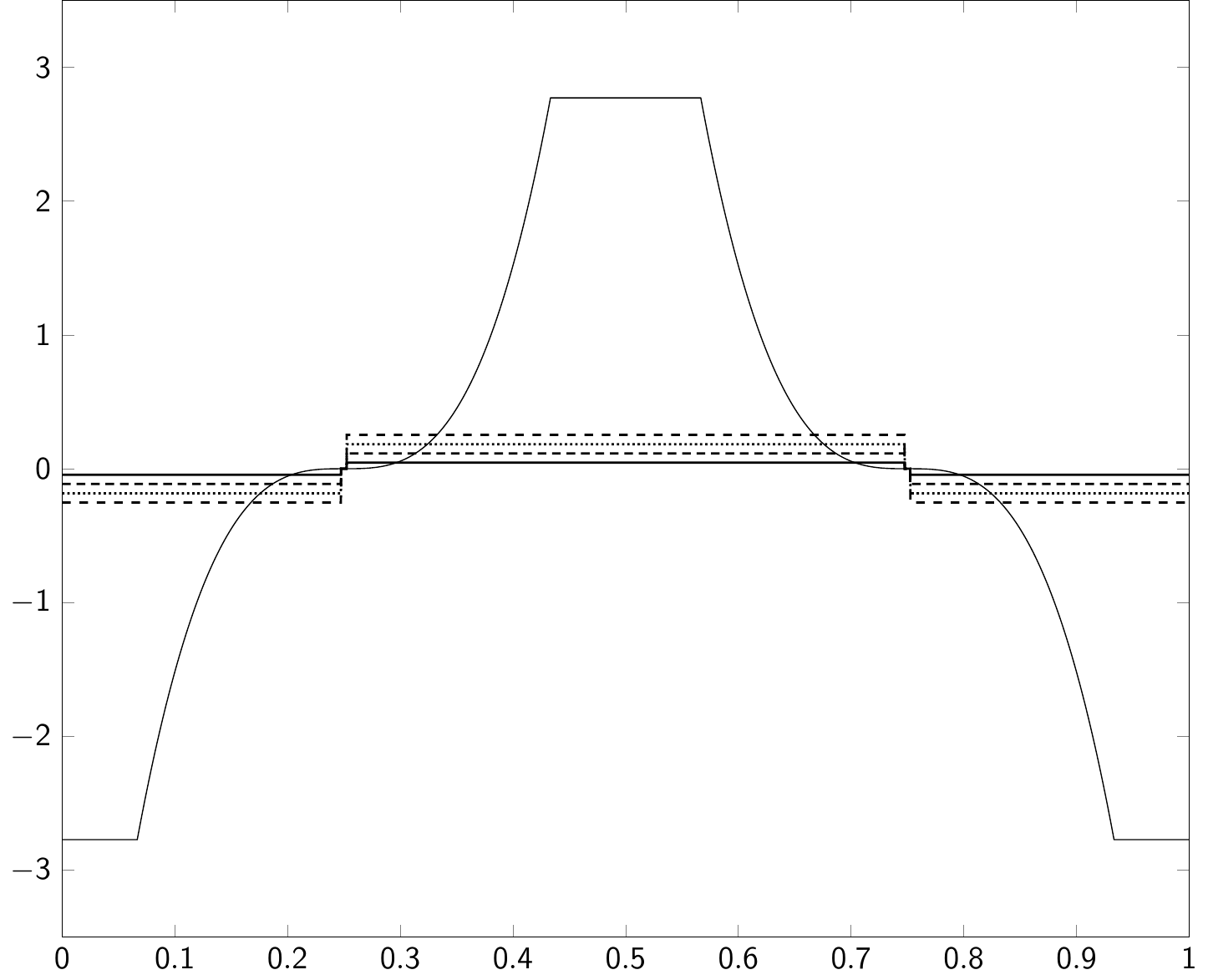}\label{Fig:grad42}}
\caption{Second numerical examples of the gradient flow.}
\label{Fig:gradient_flow2}
\end{figure}

\subsection{Example 3: Extinction time}
Our third example for \eqref{Eq:GradientFlow} is 
\begin{equation}
u^0(x) = -\cos(2\pi x),
\end{equation}
which gives
\begin{equation}
\|u^0\|_{H^{-1}_{\mathrm{av}}(\mathbb{T})} = \dfrac{1}{2\sqrt{2}\pi}.
\end{equation}
Figure \ref{Fig:gradient_flow3} shows evolution of numerical solution for third example. We use $N-200$, $\lambda=20h^{-3}$ and $\mu=30h^{-1}$ for .Figure \ref{Fig:gradient_flow3}.
Recall that our numerical scheme can compute the evolution until $\textbf{u}^k\approx \textbf{0}$. easily. Furthermore, applying the extinction time estimate \cite[Theorem 3.11]{GK11} to one-dimensional torus implies
\begin{equation}
T^*(u^0) \le C^*\|u^0\|_{H^{-1}_{\mathrm{av}}(\mathbb{T})},
\end{equation}
where $T^*(u^0)$ is the extinction time for the initial condition $u^0\in H^{-1}_{\mathrm{av}}(\mathbb{T})$ and the constant $C^*$ satisfies $\|f\|_{H^{-1}_{\mathrm{av}}(\mathbb{T})}\le C^*\int_{\mathbb{T}}|Df|$ for all $f\in H^{-1}_{\mathrm{av}}(\mathbb{T})$. It is easy to check that
\begin{equation}
\|f\|_{H^{-1}_{\mathrm{av}}(\mathbb{T})} = \left(\displaystyle\sum_{\xi\neq0}\dfrac{1}{4\pi^2}\xi^{-2}|\widehat{f}_T(\xi)|^2\right)^{1/2} \le \dfrac{1}{2\pi}\|f\|_{L^2(\mathbb{T})} \le \dfrac{1}{2\pi}\|f\|_{L^{\infty}(\mathbb{T})} \le\dfrac{1}{2\pi}\int_{\mathbb{T}}|Df|
\end{equation}
for all $f\in H^{-1}_{\mathrm{av}}(\mathbb{T})$.
Therefore, the extinction time for $u^0(x) = -\cos(2\pi x)$ can be estimated as 
\begin{equation}
T^*(u^0) \le \dfrac{1}{4\sqrt{2}\pi^2} \approx 1.7911224 \times 10^{-2}.
\end{equation}
The numerical solution is expected to be ``extinct'' in
\begin{equation}
\label{Eq:ExtinctionTime}
k \le \dfrac{T^*(u^0)}{\tau} \lessapprox  1.7911224\tau^{-1}\times 10^{-2}.
\end{equation}

\begin{table}[tb]
\centering
\begin{tabular}{|c|c|c|c|c|c|}
\hline
parameters & $\tau$ & $T^*(u^0)/\tau\lessapprox$ & $\|\textbf{u}^k\|_{\infty}< 10^{-4}$ & $\|\textbf{u}^k\|_{\infty}< 10^{-6}$ & $\|\textbf{u}^k\|_{\infty}< 10^{-8}$  \rule[0mm]{0mm}{5mm}\\\hline\hline
$N=100$, $\lambda=h^{-3}$& $10^{-6}$ &17911&4032&41769&135755 \rule[0mm]{0mm}{5mm}\\\hline
$N=100$, $\lambda=10h^{-3}$&$10^{-7}$ &179112&40311&60579&333015 \rule[0mm]{0mm}{5mm}\\\hline
$N=200$, $\lambda=10h^{-3}$&$1.25\times 10^{-8}$ &1432898&322491&592634&1267927 \rule[0mm]{0mm}{5mm}\\\hline
\end{tabular}
\caption{Time step $k$ which satisfies $\|\textbf{u}^k\|_{\infty}<10^{-4}$, $10^{-6}$ and $10^{-8}$.}
\label{Table:ExtinctionTime}
\end{table}

Table \ref{Table:ExtinctionTime} shows the time step number $k$ such that $\|\textbf{u}^k\|_{\infty}<10^{-4}$, $10^{-6}$ and $10^{-8}$ for each parameters. This result shows that we can get $\|\textbf{u}^k\|_{\infty}\lessapprox \tau$ in reasonable iteration number which is expected in \eqref{Eq:ExtinctionTime}, however, it requires more iteration to obtain smaller $\|\textbf{u}^k\|_{\infty}$.

\begin{figure}[tb]
\centering
\subfloat[][Fourth order total variation flow]{\includegraphics[clip,scale=.7]{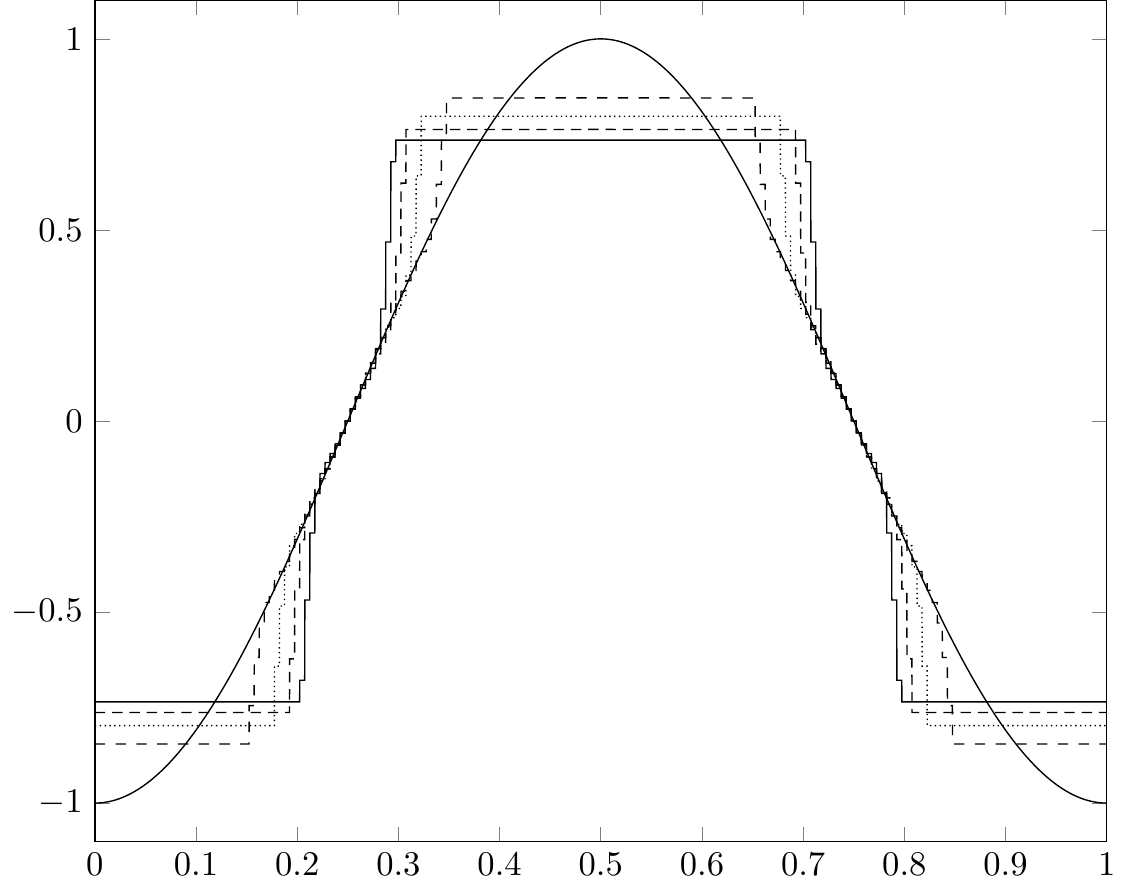}\label{Fig:gradient_flow3}}
\subfloat[][Spohn's fourth order model on $\mathbb{T}$]{\includegraphics[clip,scale=.7]{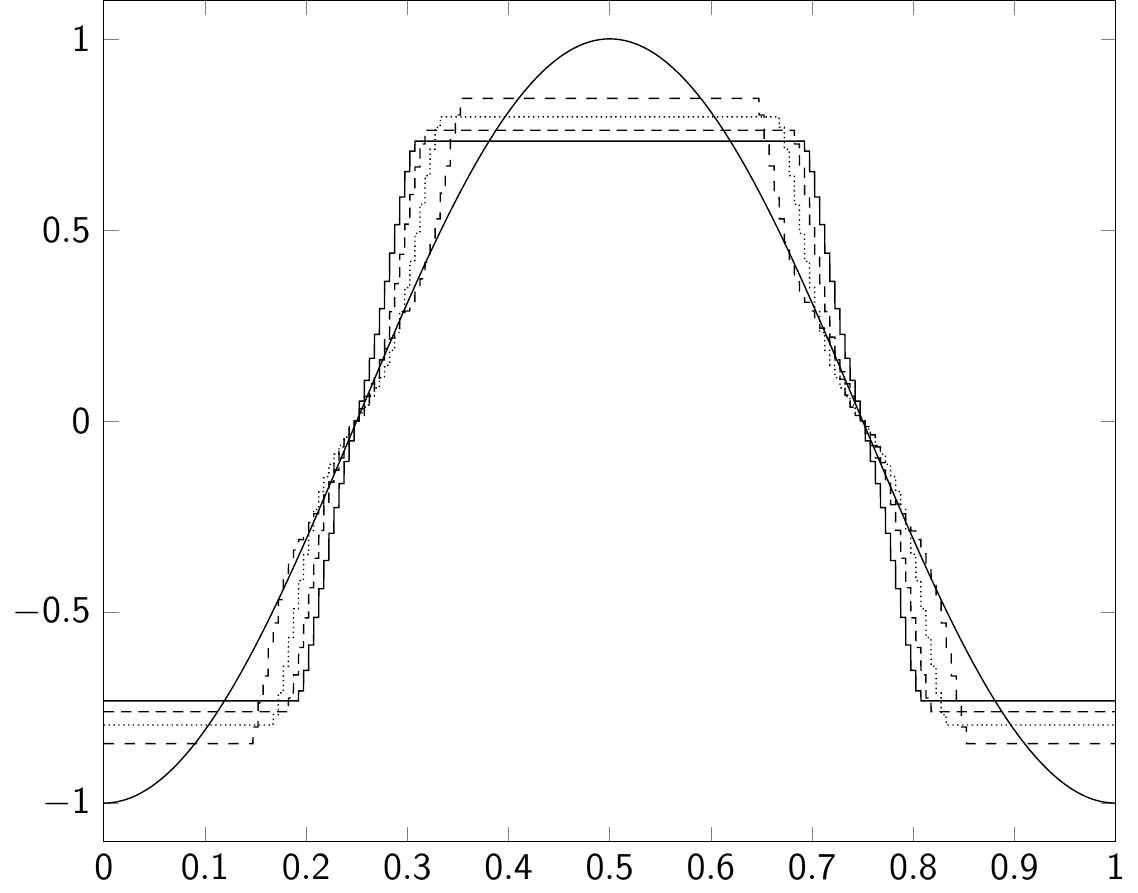}\label{Fig:gradcry}}
\caption{Numerical results for $u^0(x) = -\cos(2\pi x)$.}
\label{Fig:gradient_flow2}
\end{figure}

\subsection{Example 4: Spohn's model}
Our fourth example is split Bregman framework for Spohn's fourth order model  \eqref{Eq:SpohnModel}, which is described in Section \ref{Sec:Spohn}. 
Recall that we suppose that $p=3$ in this paper. Therefore we can apply the shrinkage operator \eqref{Eq:Shrink_for_Spohn} to split Bregman framework for Spohn's model.
Figure \ref{Fig:gradcry} shows the numerical example for $u^0(x) = -\cos(2\pi x)$, $\beta = 0.5$, $N=200$, $\lambda=50h^{-3}$ and $\mu=30h^{-1}$. 

\section{Two dimensional case}
\label{Sec:TwoDim}
The fourth order total variation flow and Spohn's model on two dimensional torus $\mathbb{T}^2$ can be computed by the similar way to one dimensional case. 
We can define $L^2_{\mathrm{av}}(\mathbb{T}^2)$, $H^1_{\mathrm{av}}(\mathbb{T}^2)$, $H^{-1}_{\mathrm{av}}(\mathbb{T}^2)$ and $(-\Delta_{\mathrm{av}})^{-1}:H^{-1}_{\mathrm{av}}(\mathbb{T}^2)\to H^1_{\mathrm{av}}(\mathbb{T}^2)$ by the generalized Fourier transform.

First, the fourth order isotropic total variation flow introduces the constraint problem
\begin{equation}
\mathop{\mathrm{minimize}}_{u\in H^{-1}_{\mathrm{av}}(\mathbb{T}^2)}\left\{\displaystyle\int_{\mathbb{T}^2}|(d_x,d_y)|+\dfrac{\lambda}{2}\|u-f\|_{H^{-1}_{\mathrm{av}}(\mathbb{T}^2)}^2 : d_x = D_xu\mbox{ and }d_y = D_yu\right\},
\end{equation}
where $D_x$, $D_y$ is distributional derivative for each variable. Note that
\begin{equation}
\|u-f\|_{H^{-1}_{\mathrm{av}}(\mathbb{T}^2)}^2=\|\nabla_x(-\Delta_{\mathrm{av}})^{-1}(u-f)\|_{L^2(\mathbb{T}^2)}^2+\|\nabla_y(-\Delta_{\mathrm{av}})^{-1}(u-f)\|_{L^2(\mathbb{T}^2)}^2,
\end{equation}
where $\nabla_x=\partial/\partial x$ and $\nabla_y = \partial/\partial y$.
Let $N_x$, $N_y$ be the partition number, $h_x=1/N_x$, $h_y = 1/N_y$, $x_n = nh_x$ and $y_n = nh_y$. Furthermore, we let $Q_{n_x,n_y} = [x_{n_x-1/2},x_{n_x+1/2})\times [y_{n_y-1/2},y_{n_y+1/2})$ and $\widehat{Q}_{n_x,n_y} = [x_{n_x-1},x_{n_x})\times [y_{n_y-1},y_{n_y})$. Then we consider the space of piecewise constant functions
\begin{subequations}
\begin{align}
V_h &= \left\{v_h :\mathbb{T}^2\to\mathbb{R} : v_h|_{Q_{n_x,n_y}}\in \mathbb{P}_0(Q_{n_x,n_y})\mbox{ for all }n_x,n_y\right\}.\\
V_{h0} &= \left\{v_h = \displaystyle\sum_{n_x=1,n_y=1}^{N_x,N_y}v_{n_x,n_y}\boldsymbol{1}_{Q_{n_x,n_y}}\in V_h : \sum_{n_x=1,n_y=1}^{N_x,N_y}v_{n_x,n_y}=0\right\},\\
\widehat{V}_h &= \left\{d_h : \Omega\to\mathbb{R} : d_h|_{\widehat{Q}_{n_x,n_y}}\in\mathbb{P}_0(\widehat{Q}_{n_x,n_y})\mbox{ for all }n_x,n_y\right\},
\end{align}
\end{subequations}
where $\Omega = [0,1)^2$. Any element $d_h\in \widehat{V}_h$ is described as $d_h = \sum_{n_x,n_y}d_{n_x,n_y}\boldsymbol{1}_{\widehat{Q}_{n_x,n_y}}$.
Let 
\begin{subequations}
\begin{align}
\widetilde{\textbf{v}} &= (v_{1,1},\dots,v_{N_x,1},v_{1,2},\dots,v_{N_x,2},\dots,v_{N_x-1,N_y}, v_{N_x,N_y})^{\mathrm{T}}\in\mathbb{R}^{N_xN_y}\\
\textbf{v} &= (v_{1,1},\dots,v_{N_x,1},v_{1,2},\dots,v_{N_x,2},\dots,v_{N_x-1,N_y})^{\mathrm{T}}\in\mathbb{R}^{N_xN_y-1}\\
\textbf{d} &= (d_{1,1},\dots,d_{N_x,1},d_{1,2},\dots,d_{N_x,2},\dots,d_{N_x,N_y})^{\mathrm{T}}\in\mathbb{R}^{N_xN_y}
\end{align}
\end{subequations}
for $v_h\in V_{h0}$ and $d_h\in \widehat{V}_h$. We define $D_{xh},D_{yh}:V_{h0}\to\widehat{V}_h\cap L^2_{\mathrm{av}}(\Omega)$ as
\begin{equation}
D_{xh}v_h = \displaystyle\sum_{n_x,n_y}(v_{n_x,n_y}-v_{n_x-1,n_y})\boldsymbol{1}_{Q_{n_x,n_y}},\quad D_{yh}v_h = \sum_{n_x,n_y}(v_{n_x,n_y}-v_{n_x,n_y-1})\boldsymbol{1}_{Q_{n_x,n_y}}.
\end{equation}
This gives
\begin{subequations}
\begin{align}
\|d_{xh}-D_{xh}u_h\|_{L^2(\Omega)}^2 &= h_xh_y\|\textbf{d}_x-h_x\nabla_{xh}R_{N_xN_y}\textbf{u}\|_2^2,\\
\|d_{yh}-D_{yh}u_h\|_{L^2(\Omega)}^2 &= h_xh_y\|\textbf{d}_y-h_y\nabla_{yh}R_{N_xN_y}\textbf{u}\|_2^2,
\end{align}
\end{subequations}
where $R_{N_xN_y}\in\mathbb{R}^{(N_xN_y)\times(N_xN_y-1)}$ is defined as equation \eqref{Def:R} and $\nabla_{xh}$, $\nabla_{yh}$ are the discrete gradient
\begin{equation}
\nabla_{xh} = h_x^{-1}I_{N_y}\otimes S_{N_x},\quad \nabla_{yh}=h_y^{-1}S_{N_y}\otimes I_{N_x},
\end{equation}
where $I_N\in\mathbb{R}^{N\times N}$ is the identity matrix and $\otimes$ is the Kronecker product.
Then our discretized problem is described as
\begin{equation*}
\begin{array}{rl}
\displaystyle\mathop{\mathrm{minimize}}_{\textbf{u}\in \mathbb{R}^{N_xN_y-1},\textbf{d}_x, \textbf{d}_y\in \mathbb{R}^{N_xN_y}}&\Biggl\{\|\textbf{d}_{xy}\|_1+\dfrac{\lambda h_xh_y}{2}\left(\|K_x(\textbf{u}-\textbf{f})\|_2^2+\|K_y(\textbf{u}-\textbf{f})\|_2^2\right)\\
&\qquad +\dfrac{\mu h_xh_y}{2}\left(\|\textbf{d}_x-h_x\nabla_{xh}R_{N_xN_y}\textbf{u}\|_2^2+\|\textbf{d}_y-h_y\nabla_{yh}R_{N_xN_y}\textbf{u}\|_2^2\right)\Biggr\},
\end{array}
\end{equation*} 
where $\textbf{d}_{xy}\in\mathbb{R}^{N_x\times N_y}$ is defined as
\begin{equation}
d_{xy,n_x,n_y} = \sqrt{d_{x,n_x,n_y}^2 + d_{y,n_x,n_y}^2}\mbox{ for all }1\le n_x\le N_x\mbox{ and }1\le n_y\le N_y,
\end{equation}
and $K_x, K_y\in\mathbb{R}^{(N_xN_y)\times(N_xN_y-1)}$ are deduced from $\nabla_x(-\Delta_{\mathrm{av}})^{-1}$ and $\nabla_y(-\Delta_{\mathrm{av}})^{-1}$, respectively. For example, we can approximate the inverse Laplacian by using
\begin{equation}
(-\Delta_{\mathrm{av}})_h = L_{N_xN_y}(\nabla_{xh}^{\mathrm{T}}\nabla_{xh}+\nabla_{yh}^{\mathrm{T}}\nabla_{yh})R_{N_xN_y}.
\end{equation}
This yields that our first scheme for two dimensional case is described as $K_x=J_x$ and $K_y=J_y$, where
\begin{equation}
J_x = \nabla_{xh}R_{N_xN_y}(-\Delta_{\mathrm{av}})_h^{-1},\quad J_y = \nabla_{yh}R_{N_xN_y}(-\Delta_{\mathrm{av}})_h^{-1}.
\end{equation}
If we let $h_x=h_y=h$, then it is required that $\lambda = O(h^{-4})$ and $\tau=O(h^{-2})$.
The split Bregman framework gives 
\begin{subequations}
\begin{empheq}[left=\empheqlbrace]{align}
&\begin{aligned}
\textbf{u}^{k+1} = \displaystyle\mathop{\mathrm{argmin}}_{\textbf{u}\in\mathbb{R}^{N_xN_y-1}}&\left\{\dfrac{\lambda h_xh_y}{2}\left(\|K_x(\textbf{u}-\textbf{f})\|_2^2+\|K_y(\textbf{u}-\textbf{f})\|_2^2\right)\right.\\
&\qquad + \dfrac{\mu h_xh_y}{2}\left(\|\textbf{d}_x^k-h_x\nabla_{xh}R_{N_xN_y}\textbf{u}-\boldsymbol{\alpha}_x^{k}\|_2^2\right.\\
&\qquad\qquad\left.\left.+\|\textbf{d}_y^k-h_y\nabla_{yh}R_{N_xN_y}\textbf{u}-\boldsymbol{\alpha}_y^{k}\|_2^2\right)\right\},
\end{aligned}\label{Eq:2Du1}\\
&\begin{aligned}
(\textbf{d}_x^{k+1},\textbf{d}_y^{k+1}) = \displaystyle\mathop{\mathrm{argmin}}_{\textbf{d}_x,\textbf{d}_y\in\mathbb{R}^{N_xN_y}}&\left\{\|\textbf{d}_{xy}\|_1+\dfrac{\mu h_xh_y}{2}\left(\|\textbf{d}_x-h_x\nabla_{xh}R_{N_xN_y}\textbf{u}^{k+1}-\boldsymbol{\alpha}_x^{k}\|_2^2\right.\right.\\
&\left.\left.\qquad+\|\textbf{d}_y-h_y\nabla_{yh}R_{N_xN_y}\textbf{u}^{k+1}-\boldsymbol{\alpha}_y^{k}\|_2^2\right)\right\},
\end{aligned}\label{Eq:2Dd1}\\
&\boldsymbol{\alpha}_x^{k+1} = \boldsymbol{\alpha}_x^k-\textbf{d}_x^{k+1}+h_x\nabla_{xh}R_{N_xN_y}\textbf{u}^{k+1},\quad\boldsymbol{\alpha}_y^{k+1} = \boldsymbol{\alpha}_y^k-\textbf{d}_y^{k+1}+h_y\nabla_{yh}R_{N_xN_y}\textbf{u}^{k+1},
\end{empheq}
\end{subequations}
where $\textbf{f}\in\mathbb{R}^{N_xN_y-1}$ is given data or $\textbf{f}=\textbf{u}^k$, $\boldsymbol{\alpha}_x^0=\boldsymbol{\alpha}_x^0 = \textbf{0}$, $\textbf{u}^0$ is given as $\textbf{0}$ or initial condition, and $\textbf{d}_x^0 = h_x\nabla_{xh}R_{N_xN_y}\textbf{u}^0$, $\textbf{d}_y^0 = h_y\nabla_{yh}R_{N_xN_y}\textbf{u}^0$.

Note that the equation \eqref{Eq:2Dd1} is essentially the same formulation as the one of split Bregman framework for second order isotropic problem, which is mentioned in \cite{GO09}. 
The Euler-Lagrange equation for equation \eqref{Eq:2Dd1} yields
\begin{subequations}
\begin{align}
\dfrac{(\textbf{d}_x^{k+1})_n}{|(\textbf{d}_{xy}^{k+1})_n|}+\mu h_xh_y\left(\textbf{d}_x^{k+1}-h_x\nabla_{xh}R_{N_xN_y}\textbf{u}^{k+1}-\boldsymbol{\alpha}_x^k\right)_n=0,\label{Eq:Isotropic_ite_dx}\\
\dfrac{(\textbf{d}_y^{k+1})_n}{|(\textbf{d}_{xy}^{k+1})_n|}+\mu h_xh_y\left(\textbf{d}_y^{k+1}-h_y\nabla_{yh}R_{N_xN_y}\textbf{u}^{k+1}-\boldsymbol{\alpha}_y^k\right)_n=0\label{Eq:Isotropic_ite_dy}
\end{align}
\end{subequations}
for all $n=1,\dots, N_xN_y$.
We consider the approximation
\begin{equation}
\label{Eq:2D_singular_approx}
\dfrac{(\textbf{d}_x^{k+1})_n}{|(\textbf{d}_{xy}^{k+1})_n|} \approx \dfrac{(\textbf{d}_x^{k+1})_n}{|(\textbf{d}_{x}^{k+1})_n|}\cdot\dfrac{|s_{x,n}^k|}{s_n^k},\quad\dfrac{(\textbf{d}_y^{k+1})_n}{|(\textbf{d}_{xy}^{k+1})_n|} \approx \dfrac{(\textbf{d}_y^{k+1})_n}{|(\textbf{d}_{y}^{k+1})_n|}\cdot\dfrac{|s_{y,n}^k|}{s_n^k},
\end{equation}
where
\begin{equation*}
s_n^k = \sqrt{(s_{x,n}^k)^2+(s_{y,n}^k)^2},\quad s_{x,n}^k = (h_x\nabla_{xh}R_{N_xN_y}\textbf{u}^{k+1}+\boldsymbol{\alpha}_x^k)_n,\quad s_{y,n}^k = (h_y\nabla_{yh}R_{N_xN_y}\textbf{u}^{k+1}+\boldsymbol{\alpha}_y^k)_n.
\end{equation*}
Applying them into equations \eqref{Eq:Isotropic_ite_dx} and \eqref{Eq:Isotropic_ite_dy} give the following shrinkage formula, which are equivalent to ones of \cite[Section 4.1]{GO09}:
\begin{equation*}
(d_x^{k+1})_n = \dfrac{s_{x,n}^k}{|s_{x,n}^k|}\max\left\{|s_{x,n}^k|-\dfrac{|s_{x,n}^k|}{\mu h_xh_ys_n^k},0\right\},\quad (d_y^{k+1})_n = \dfrac{s_{y,n}^k}{|s_{y,n}^k|}\max\left\{|s_{y,n}^k|-\dfrac{|s_{y,n}^k|}{\mu h_xh_ys_n^k},0\right\}.
\end{equation*}
Figure \ref{Fig:2D_Isotropic} shows the numerical result of fourth order isotropic total variation flow \eqref{Eq:FourthOrderTVFlow} in $\mathbb{T}^2$ with initial data $u^0(x,y) = x(x-1)y(y-1)-1/36$. We use $N_x=N_y=40$, $\lambda=5h^{-4}$ and $\mu=20h^{-2}$.

Next, the fourth order anisotropic total variation flow 
\begin{equation}
u_t = -\Delta\left(\operatorname{div}\left(\dfrac{\nabla_xu}{|\nabla_xu|},\dfrac{\nabla_yu}{|\nabla_yu|}\right)\right).
\end{equation}
\begin{figure}[H]
\centering
\subfloat[][Fourth order isotropic total variation flow]{\includegraphics[clip,scale=1.5]{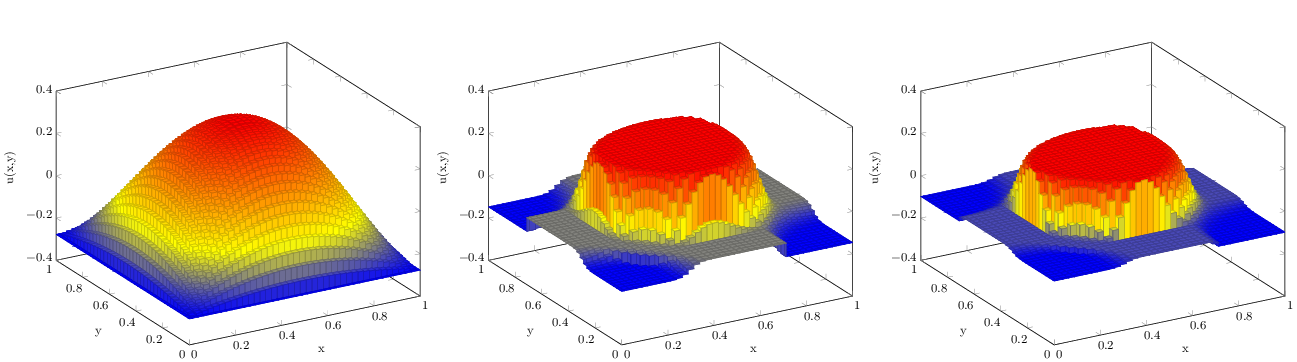}\label{Fig:2D_Isotropic}}
\\
\subfloat[][Fourth order anisotropic total variation flow]{\includegraphics[clip,scale=1.5]{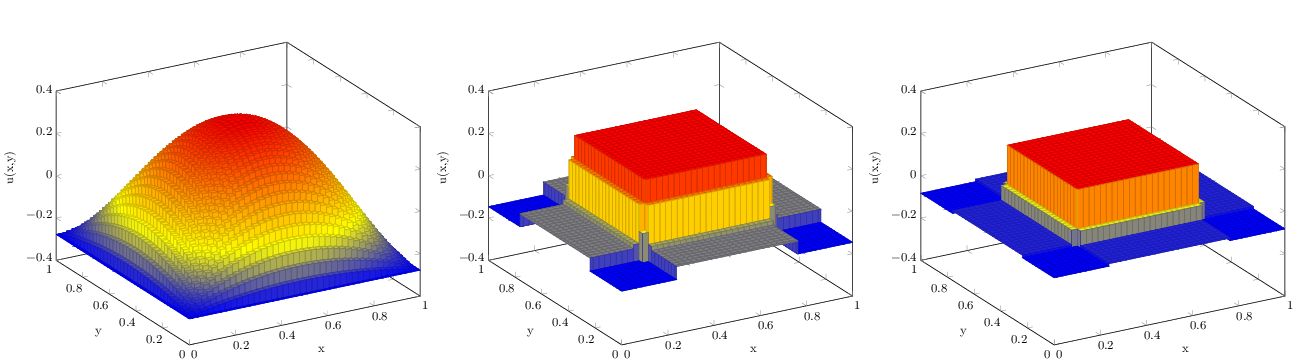}\label{Fig:2D_Anisotropic}}
\\
\subfloat[][Spohn's fourth order model on $\mathbb{T}^2$]{\includegraphics[clip,scale=1.5]{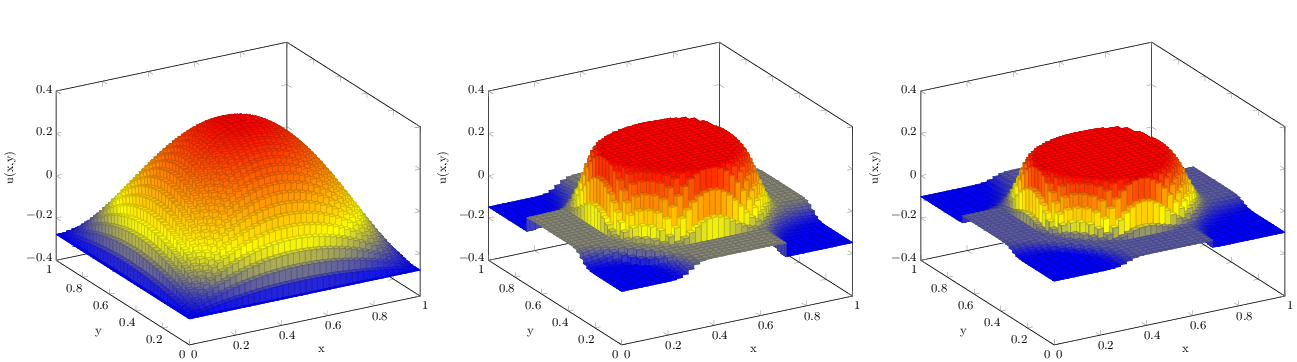}\label{Fig:2D_Spohn}}
\caption{Numerical results of two-dimensional problems.}
\end{figure}
\noindent 
Letting $F(u) = \int_{\mathbb{T}^2}\left(|D_xu|+|D_yu|\right)$ implies that formally we have
\begin{align*}
\left(\Delta\left(\operatorname{div}\left(\dfrac{\nabla_xu}{|\nabla_xu|},\dfrac{\nabla_yu}{|\nabla_yu|}\right)\right),v-u\right)_{H^{-1}_{\mathrm{av}}(\mathbb{T}^2)} &= \left(-\operatorname{div}\left(\dfrac{\nabla_xu}{|\nabla_xu|},\dfrac{\nabla_yu}{|\nabla_yu|}\right),v-u\right)_{L^2_{\mathrm{av}}(\mathbb{T}^2)}\\
&= \displaystyle\int_{\mathbb{T}^2}\left(\dfrac{\nabla_xu\overline{\nabla_xv}}{|\nabla_xu|}-|\nabla_xu|+\dfrac{\nabla_yu\overline{\nabla_yv}}{|\nabla_yu|}-|\nabla_yu|\right)\\
&\le F(v)-F(u),
\end{align*}
therefore $u_t\in-\partial_{H^{-1}_{\mathrm{av}}(\mathbb{T}^2)} F$. We apply the backward Euler method and obtain
\begin{equation}
u^{k+1} = \mathop{\mathrm{argmin}}_{u\in H^{-1}_{\mathrm{av}}(\mathbb{T}^2)}\left\{\displaystyle\int_{\mathbb{T}^2}\left(|D_xu|+|D_yu|\right)+\dfrac{1}{2\tau}\|u-u^k\|_{H^{-1}_{\mathrm{av}}(\mathbb{T}^2)}^2\right\},
\end{equation}
which introduces the constraint problem
\begin{equation}
\mathop{\mathrm{minimize}}_{u\in H^{-1}_{\mathrm{av}}(\mathbb{T}^2)}\left\{\displaystyle\int_{\mathbb{T}^2}\left(|d_x|+|d_y|\right)+\dfrac{\lambda}{2}\|u-f\|_{H^{-1}_{\mathrm{av}}(\mathbb{T}^2)}^2 : d_x = D_xu\mbox{ and }d_y = D_yu\right\},
\end{equation}
This, combining with the split Bregman framework gives
\begin{subequations}
\begin{empheq}[left=\empheqlbrace]{align}
&\begin{aligned}
\textbf{u}^{k+1} = \displaystyle\mathop{\mathrm{argmin}}_{\textbf{u}\in\mathbb{R}^{N_xN_y-1}}&\left\{\dfrac{\lambda h_xh_y}{2}\left(\|K_x(\textbf{u}-\textbf{f})\|_2^2+\|K_y(\textbf{u}-\textbf{f})\|_2^2\right)\right.\\
&\qquad + \dfrac{\mu h_xh_y}{2}\left(\|\textbf{d}_x^k-h_x\nabla_{xh}R_{N_xN_y}\textbf{u}-\boldsymbol{\alpha}_x^{k}\|_2^2\right.\\
&\qquad\qquad\left.\left.+\|\textbf{d}_y^k-h_y\nabla_{yh}R_{N_xN_y}\textbf{u}-\boldsymbol{\alpha}_y^{k}\|_2^2\right)\right\},
\end{aligned}\label{Eq:2Du2}\\
&\textbf{d}_x^{k+1} = \displaystyle\mathop{\mathrm{argmin}}_{\textbf{d}_x\in\mathbb{R}^{N_xN_y}}\left\{\|\textbf{d}_x\|_1+\dfrac{\mu h_xh_y}{2}\|\textbf{d}_x-h_x\nabla_{xh}R_{N_xN_y}\textbf{u}^{k+1}-\boldsymbol{\alpha}_x^{k}\|_2^2\right\}\label{Eq:2D_anisotropic_dx},\\
&\textbf{d}_y^{k+1} = \displaystyle\mathop{\mathrm{argmin}}_{\textbf{d}_y\in\mathbb{R}^{N_xN_y}}\left\{\|\textbf{d}_y\|_1+\dfrac{\mu h_xh_y}{2}\|\textbf{d}_y-h_y\nabla_{yh}R_{N_xN_y}\textbf{u}^{k+1}-\boldsymbol{\alpha}_y^{k}\|_2^2\right\}\label{Eq:2D_anisotropic_dy},\\
&\boldsymbol{\alpha}_x^{k+1} = \boldsymbol{\alpha}_x^k-\textbf{d}_x^{k+1}+h_x\nabla_{xh}R_{N_xN_y}\textbf{u}^{k+1},\quad\boldsymbol{\alpha}_y^{k+1} = \boldsymbol{\alpha}_y^k-\textbf{d}_y^{k+1}+h_y\nabla_{yh}R_{N_xN_y}\textbf{u}^{k+1}.
\end{empheq}
\end{subequations}
We can apply the shrinking method \eqref{Op:shrink} to equations \eqref{Eq:2D_anisotropic_dx} and \eqref{Eq:2D_anisotropic_dy}.
Figure \ref{Fig:2D_Anisotropic} presents the evolution of fourth order anisotropic total variation flow for $u^0(x,y) = x(x-1)y(y-1)-1/36$, $N_x=N_y=40$, $\lambda=5h^{-4}$ and $\mu=20h^{-2}$. 
For second order anisotropic total variation flow, \L asica, Moll and Mucha \cite{LMM17} have considered rectangular domain $\Omega\subset\mathbb{R}^2$ or $\Omega=\mathbb{R}^2$ and rigorously proved that if the initial profile is piecewise constant, then the exact solution is piecewise constant. 
We can infer from our numerical experiment \ref{Fig:2D_Anisotropic} that their theoretical result is true also for fourth order anisotropic total variation flow.

Finally, we consider two dimensional Spohn's fourth order model. The split Bregman framework provides
\begin{subequations}
\begin{empheq}[left=\empheqlbrace]{align}
&\begin{aligned}
\textbf{u}^{k+1} = \displaystyle\mathop{\mathrm{argmin}}_{\textbf{u}\in\mathbb{R}^{N_xN_y-1}}&\left\{\dfrac{\lambda h_xh_y}{2}\left(\|K_x(\textbf{u}-\textbf{f})\|_2^2+\|K_y(\textbf{u}-\textbf{f})\|_2^2\right)\right.\\
&\qquad + \dfrac{\mu h_xh_y}{2}\left(\|\textbf{d}_x^k-h_x\nabla_{xh}R_{N_xN_y}\textbf{u}-\boldsymbol{\alpha}_x^{k}\|_2^2\right.\\
&\qquad\qquad\left.\left.+\|\textbf{d}_y^k-h_y\nabla_{yh}R_{N_xN_y}\textbf{u}-\boldsymbol{\alpha}_y^{k}\|_2^2\right)\right\},
\end{aligned}\\
&\begin{aligned}
(\textbf{d}_x^{k+1},\textbf{d}_y^{k+1}) = \displaystyle\mathop{\mathrm{argmin}}_{\textbf{d}_x,\textbf{d}_y\in\mathbb{R}^{N_xN_y}}&\left\{\beta\|\textbf{d}_{xy}\|_1+\dfrac{1}{p}\|\textbf{d}_{xy}\|_p^p\right.\\
&\qquad+\dfrac{\mu h_xh_y}{2}\left(\|\textbf{d}_x-h_x\nabla_{xh}R_{N_xN_y}\textbf{u}^{k+1}-\boldsymbol{\alpha}_x^{k}\|_2^2\right.\\
&\qquad\qquad\left.\left.+\|\textbf{d}_y-h_y\nabla_{yh}R_{N_xN_y}\textbf{u}^{k+1}-\boldsymbol{\alpha}_y^{k}\|_2^2\right)\right\},
\end{aligned}\label{Eq:2DSpohn_d}\\
&\boldsymbol{\alpha}_x^{k+1} = \boldsymbol{\alpha}_x^k-\textbf{d}_x^{k+1}+h_x\nabla_{xh}R_{N_xN_y}\textbf{u}^{k+1},\quad\boldsymbol{\alpha}_y^{k+1} = \boldsymbol{\alpha}_y^k-\textbf{d}_y^{k+1}+h_y\nabla_{yh}R_{N_xN_y}\textbf{u}^{k+1}.
\end{empheq}
\end{subequations}
The Euler-Lagrange equation for \eqref{Eq:2DSpohn_d} can be approximated by equation \eqref{Eq:2D_singular_approx}.
In this paper, we always suppose that $p=3$. Note that the approximation \eqref{Eq:2D_singular_approx} implies
\begin{equation}
|(\textbf{d}_{xy}^{k+1})_n| \approx |(\textbf{d}_x^{k+1})_n|\cdot\dfrac{s_n^k}{|s_{x,n}^k|}\quad\mbox{ and }\quad|(\textbf{d}_{xy}^{k+1})_n| \approx |(\textbf{d}_y^{k+1})_n|\cdot\dfrac{s_n^k}{|s_{y,n}^k|}.
\end{equation}
We obtain approximated Euler-Lagrange equations
\begin{subequations}
\begin{align}
\beta\dfrac{(\textbf{d}_x^{k+1})_n}{|(\textbf{d}_{x}^{k+1})_n|}\cdot\dfrac{|s_{x,n}^k|}{s_n^k}+(\textbf{d}_x^{k+1})_n|(\textbf{d}_x^{k+1})_n|\cdot\dfrac{s_n^k}{|s_{x,n}^k|}+\mu h_xh_y((\textbf{d}_x^{k+1})_n-s_{x,n}^k)=0,\label{Eq:Spohn_ite_dx}\\
\beta\dfrac{(\textbf{d}_y^{k+1})_n}{|(\textbf{d}_{y}^{k+1})_n|}\cdot\dfrac{|s_{y,n}^k|}{s_n^k}+(\textbf{d}_y^{k+1})_n |(\textbf{d}_y^{k+1})_n|\cdot\dfrac{s_n^k}{|s_{y,n}^k|}+\mu h_xh_y((\textbf{d}_y^{k+1})_n-s_{y,n}^k)=0\label{Eq:Spohn_ite_dy}.
\end{align}
\end{subequations}

By the similar way to one dimensional case, we provide the shrinkage operators of the form
\begin{subequations}
\begin{align}
(\textbf{d}_x^{k+1})_n &= \dfrac{\mu h_xh_y|s_{x,n}^k|}{2s_n^k}\cdot\dfrac{s_{x,n}^k}{|s_{x,n}^k|}\left(-1+\sqrt{1+\dfrac{4s_n^k}{\mu h_xh_y|s_{x,n}^k|}\max\left\{|s_{x,n}^k|-\dfrac{\beta |s_{x,n}^k|}{\mu h_xh_ys_n^k},0\right\}}\right),\\
(\textbf{d}_y^{k+1})_n &= \dfrac{\mu h_xh_y|s_{y,n}^k|}{2s_n^k}\cdot\dfrac{s_{y,n}^k}{|s_{y,n}^k|}\left(-1+\sqrt{1+\dfrac{4s_n^k}{\mu h_xh_y|s_{y,n}^k|}\max\left\{|s_{y,n}^k|-\dfrac{\beta |s_{y,n}^k|}{\mu h_xh_ys_n^k},0\right\}}\right).
\end{align}
\end{subequations}
Figure \ref{Fig:2D_Spohn} shows the numerical result of split Bregman framework for Spohn's forth order model. We use $p=3$, $\beta=0.25$, $N_x=N_y=40$, $\lambda=1.25h^{-4}$ and $\mu=5h^{-2}$. Moreover, we use the initial value $u^0(x,y) = x(x-1)y(y-1)-1/36$, which is considered in \cite{KV10}. We can obtain the similar numerical result quite effectively by split Bregman framework.

%
%


\section{Conclusion}
In this study, we propose a new numerical scheme for the OSV model, fourth order total variation flow and Spohn's fourth order model. Our scheme is based on the split Bregman framework for the ROF model and second order total variation flow. We demonstrate several numerical examples for one dimensional and two dimensional problems under periodic boundary condition. We use the parameters $\lambda=O(h^{-3})$, $\mu=O(h^{-1})$ for one dimensional case, and $\lambda=O(h^{-4})$, $\mu=O(h^{-2})$ for two dimensional case. For fourth order total variation flow, our numerical results approximately represent the flat facet and discontinuity, which is expected by the theoretical result for the exact profile. Furthermore, we propose new shrinkage operators for Spohn's model. Numerical results for Spohn's model show facet and relaxation.
 
\section*{Acknowledgement}
A part of the work of the second author was done when he was a postdoc fellow at the University of Tokyo. Its hospitality is gratefully acknowledged. The work of the first author was partly supported by the Japan Society for the Promotion of Science through the grant No. 26220702 (Kiban S), No. 19H00639 (Kiban A), No. 18H05323 (Kaitaku), No. 17H01091 (Kiban A) and No. 16H03948 (Kiban B).

\bibliographystyle{plain}

\end{document}